\documentclass[AMS,STIX1COL]{WileyNJD-v2}

\usepackage{bm}
\usepackage{itm_symbol}
\usepackage{xfrac}
\usepackage{booktabs}
\usepackage{makecell}
\usepackage{siunitx}
\usepackage[utf8]{inputenc}
\articletype{Review}%

\received{20 September 2022} %
\revised{6 March 2023} %
\accepted{23 March 2023} %

\newcommand{\mkr}{\mathbf{x}}

\newcommand{\fopt}{f^{\mathrm{M}}_{\circ}}

\newcommand{\fmN}{\tilde{f}^{\mathrm{M}}}
\newcommand{\lmN}{\tilde{\ell}^{\mathrm{M}}}
\newcommand{\lopt}{\ell^{\mathrm{M}}_{\circ}}
\newcommand{\ltslk}{\ell^{\mathrm{T}}_{s}}
\newcommand{\pen}{\alpha}

\newcommand{\tcA}{\tau_{\mathrm{A}}}
\newcommand{\tcD}{\tau_{\mathrm{D}}}

\newcommand{\vmN}{\tilde{v}^{\mathrm{M}}}
\newcommand{\vmax}{v^{\mathrm{M}}_{\circ}}
\newcommand{\fvN}{\mathbf{\mathrm{f}^{\mathrm{V}}}(\tilde{v}^{\mathrm{M}}) }
\newcommand{\flN}{ \mathbf{ \mathrm{f}^{\mathrm{L}} } (\tilde{\ell}^{\mathrm{M}}) }
\newcommand{\fpeN}{ \mathbf{ \mathrm{f}^{\mathrm{PE}} } (\tilde{\ell}^{\mathrm{M}}) }
\newcommand{\ftN}{ \mathbf{ \mathrm{f}^{\mathrm{T}} } (\tilde{\ell}^{\mathrm{T}}) }
\newcommand{\act}{a}
\newcommand{\dact}{\dot{a}}

\newcommand{\vmkr}{\mathbf{z}}

\newcommand{\meh}[1]{\mathrm{#1}}
\newcommand\thefont{\expandafter\string\the\font}
\hypersetup{
	hidelinks, %
	pdftitle={Dynamic Human Body Models in Vehicle Safety: An Overview} %
}
\newcolumntype{L}[1]{>{\raggedright\let\newline\\\arraybackslash\hspace{0pt}}p{#1}}

\begin{document}

\title{Dynamic Human Body Models in Vehicle Safety: An Overview} %

\author[1]{N. Fahse}
\author[1,3]{M. Millard}
\author[1]{F. Kempter}
\author[1]{S. Maier}
\author[2]{M. Roller}
\author[1]{J. Fehr*}

\authormark{Fahse \textsc{et al}}

\address[1]{\orgdiv{Institute of Engineering and Computational Mechanics (ITM)}, \orgname{University of Stuttgart}, \orgaddress{\state{Stuttgart}, \country{Germany}}}

\address[2]{\orgdiv{Fraunhofer Institute for Industrial Mathematics (ITWM)}, \orgname{Fraunhofer-Gesellschaft}, \orgaddress{\state{Kaiserslautern}, \country{Germany}}}

\address[3]{\orgdiv{Institute of Sport and Movement Science (Inspo)}, \orgname{University of Stuttgart}, \orgaddress{\state{Stuttgart}, \country{Germany}}}

\corres{*J. Fehr, Institute of Engineering and Computational Mechanics (ITM), University of Stuttgart, Stuttgart, Germany.\\\email{joerg.fehr@itm.uni-stuttgart.de}\\\newline
N.Fahse and M.Millard are equal first authors}

\abstract[Summary]{
	Significant trends in the vehicle industry are autonomous driving, micromobility, electrification and the increased use of shared mobility solutions.
	These new vehicle automation and mobility classes lead to a larger number of occupant positions, interiors and load directions.
	As safety systems interact with and protect occupants, it is essential to place the human, with its variability and vulnerability, at the center of the design and operation of these systems.
	Digital human body models (HBMs) can help meet these requirements and are therefore increasingly being integrated into the development of new vehicle models.
	This contribution provides an overview of current HBMs and their applications in vehicle safety in different driving modes.
	The authors briefly introduce the underlying mathematical methods and present a selection of HBMs to the reader.
	An overview table with guideline values for simulation times, common applications and available variants of the models is provided.
	To provide insight into the broad application of HBMs, the authors present three case studies in the field of vehicle safety:
	(i)~in-crash finite element simulations and injuries of riders on a motorcycle;
	(ii)~scenario-based assessment of the active pre-crash behavior of occupants with the Madymo multibody HBM;
	(iii)~prediction of human behavior in a take-over scenario using the EMMA model.
}


\jnlcitation{\cname{%
\author{N. Fahse}, 
\author{M. Millard},
\author{F. Kempter}, 
\author{S. Maier}, 
\author{M. Roller}, 
 and 
\author{J. Fehr}} (\cyear{2022}), 
\ctitle{Dynamic Human Body Models in Vehicle Safety: An Overview}, \cjournal{TODO}, \cvol{TODO}.}

\maketitle

\section{Introduction\label{sec:Introduction}}

Mary Ward was the first victim of a traffic accident when she was killed in 1869 by a steam-powered automobile~\cite{FallonONeill05}. 
Today vehicles are much safer~\cite{NSC22} due to improvements in traffic regulations, the structural design of cars, and the addition of specific safety systems. 
The structural design of a modern vehicle is heavily influenced by safety concerns. 
For example, the chassis is designed to crumple and absorb energy during an accident, the doors have been stiffened to prevent passengers from being crushed during an impact from the side, and the hood is designed to flex and cushion the impact of a pedestrian. 
In addition, modern vehicles include several safety systems that did not exist in Mary Ward's time: pre-tensioned seat belts, front airbags, and side airbags.

The development of many vehicle safety systems has been supported by experimental and computerized crash simulations that include representations of the human form.
During an experimental crash test, an anthropomorphic test device (ATD) is used to predict the response of a human passenger. 
ATDs are instrumented passive mechanical devices that can mimic the geometry, mass distribution, joints and passive joint stiffness of the human body.
These test devices have been used to develop many safety systems, and remain part of the certification process of modern vehicles~\cite{EuroNCAP21, NHTSA22}.
Unfortunately, physical crash tests with ATDs are expensive, cannot be done without a physical prototype of the car, and can only provide an indication of an injury.
In contrast, computer models of the human body can be simulated without a physical prototype of a car and capture enough anatomical detail to make specific tissue-level injury predictions.
The explosion in the number of different human models and the variety of simulation methods has motivated us to write this review to help people working in the field make informed decisions when choosing both a model and a method to answer a specific vehicle safety problem.

There are a wide variety of different types of digital human body models (HBMs) \footnote{In the remainder of the paper the term HBM will be used to refer to a mathematical model of the human body that is simulated using a computer.} that are needed to address the complexities of car design: there are ergonomic models used to assess comfort, rigid-body models that are used to predict movement while driving, and detailed musculoskeletal and finite-element models used to predict injury during car accidents. 
While ATDs must be physically robust to produce reproducible results and are uniaxial by design, HBMs can use anatomically accurate joints and tissue models.
In contrast, HMBs serve to act as a virtual microscope to study reachability, kinematics, comfort, and injury mechanisms.
The number and variety of models is likely to increase in response to the greater variety of vehicles that are being developed, and the increasing number of passenger types that must be considered. 

\begin{figure}[t]
	\centerline{\includegraphics[scale=1]{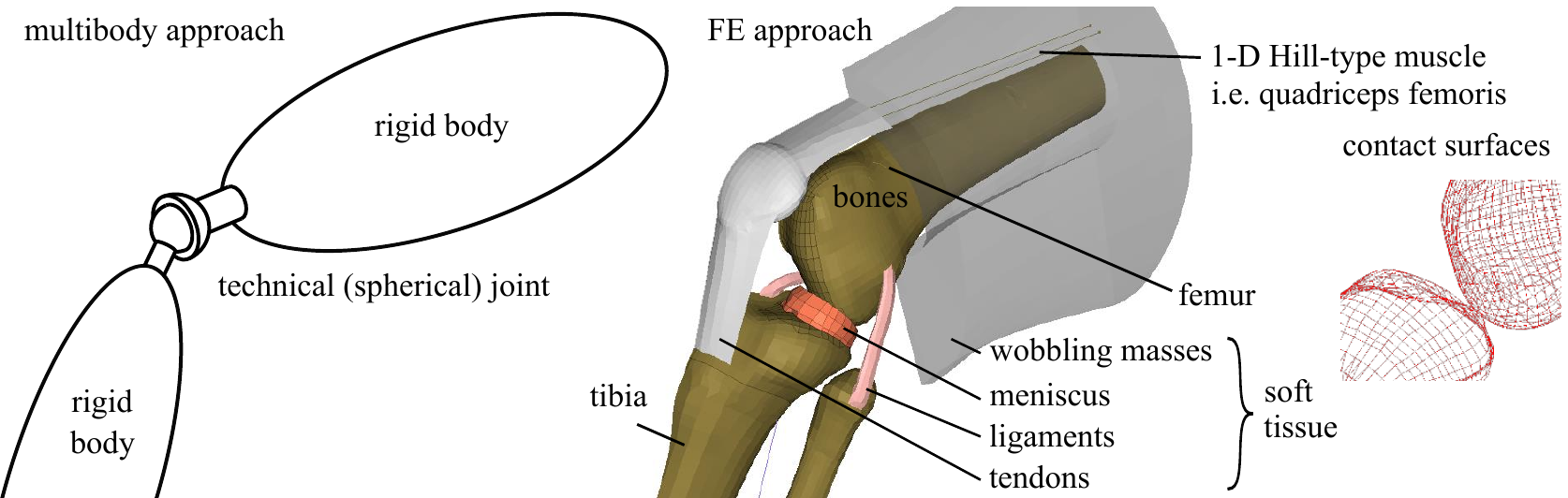}}
	\caption{\enspace On the left a simple multibody system serves as an approximation of the knee joint. On the right in a detailed FE knee model approach, the femur and tibia bones may be connected by the ligaments and tendons of the knee, separated by the meniscus and other articular cartilage, driven by the muscles of the leg, and surrounded by wobbling masses elements. Note that only a few of the ligaments of the knee are depicted and the articular capsule is not visible.\label{fig:Knee_Joint}}	
\end{figure}

Over the past decade, several HBMs with an increasing level of detail were developed to predict human response and assess tissue-level injury risk~\cite{BullingerHoffmannMuhlstedt16, ScatagliniPaul19}.
For a first approximation, the human body can be modeled as a rigid multibody system (MBS) that gets its mobility from idealized joints that connect each segment (an example featuring the knee joint appears in the left panel Fig.~\ref{fig:Knee_Joint}) as is done in the model SAMMIE (System for Aiding Man Machine Interaction Evaluation)~\cite{CaseEtAl16, CaseEtAl86}.
Additional detail can be added to a rigid body model by replacing idealized joints with anatomically accurate representations that include ligaments, tendons, cartilage, and the network of muscles that surround a joint (right panel Fig.~\ref{fig:Knee_Joint}).
Rigid body models, however, have limitations: bones are assumed to be rigid and tissues are spatially lumped.
In contrast, finite element (FE) models can be used to accurately capture both the geometry and material properties of a variety of tissues making it possible to predict injury~\cite{SanderEtAl13}.

Typically, MBS models have fewer degrees of freedom (DOFs) than FE models.
Within the context of an individual MBS model, different coordinate systems can be used to make equivalent versions of the same model.
In a redundant coordinate system each body segment has three linear and three rotational coordinates and joints are represented mathematically as constraints and have equations of motion that are differential algebraic equations (DAEs) in structure. 
While DAEs can be simulated forward in time, applying inverse kinematics and inverse dynamics methods to a system of DAEs is challenging.
As a benefit, however, redundant coordinate systems make it possible to use anatomically detailed joint models that include ligaments and compliant cartilage.
In contrast, a minimal coordinate system treats the joint coordinates as the DOFs and results in a smaller set of equations.
In the special case that the model has a tree-like structure and contains no closed-loops, the equations of motion take the form of a system of ordinary differential equations (ODEs) that can be easily used with both forward and inverse methods.

In either case, the coordinate system ultimately does not affect the $f$ DOFs of the model: specifying $f$ generalized (minimal) coordinates $\yb$ of a model will uniquely specify the position of the model in space.
The equations of motion for a holonomic MBS with a tree-like structure are
\begin{equation}
\label{Eq:Mubo_equation}
\boldsymbol{M}(\yb,t)\cdot\ybpp={\boldsymbol{k}}_{\mathrm{c}} (\yb,\ybp,t)+ {\boldsymbol{k}}_{\mathrm{e}} (\yb,\ybp,t), 
\end{equation}
where $\boldsymbol{M} \in \mathbb{R}^{f\times f}$ represents the global mass matrix, ${\boldsymbol{k}}_{\mathrm{c}} \in \mathbb{R}^{f\times 1}$ the vector of the generalized Coriolis and centrifugal forces and ${\boldsymbol{k}}_{\mathrm{e}} \in \mathbb{R}^{f\times 1}$ the vector of the generalized applied forces~\cite{Woernle11}. 

While MBS models are accurate enough for simulating voluntary motion, FE models are needed to accurately simulate the movements and injuries that take place during a car crash.
The human body can be described as an elasto-plastic dynamic continuum system and modeled using the FE method~\cite{BelytschkoLiuMoran00}. 
The FE approach allows the representation of nonlinear material behavior, large deformations, and multiphysics coupling.
The FE equation 
\begin{equation}\label{eq:Nonlinear_Fe_Equation}
\Mb\cdot\qbpp+\fb_\mathrm{int}(\qb,\qbp)=\fb_\mathrm{ext}, 
\end{equation}
consists of the symmetric and usually diagonal mass matrix $\Mb\in \mathbb{R}^{N\times N}$;
the nodal coordinates $\qb$, nodal velocities $\qbp$, and nodal accelerations $\qbpp$.
The vector of nodal forces from the internal forces (element stresses) comes from the material laws (the relation between stress, strain, strain rate and loading history) $\fb_\mathrm{int}(\qb,\qbp)$
and the vector of the external surface and volume forces $\fb_\mathrm{ext}\in \mathbb{R}^{N}$. 
The number of DOFs, $N$, of an FE model is given by the product of the number of nodes in the mesh with the number of DOFs of each node. The number of DOFs per node can vary widely depending on the type of element being used.
States associated with active tissues, such as muscle, further increase the DOF count of the model.

Active muscle models are being included in HBMs (creating an active HBM, or AHBM) due to the potential influence of muscle activity on the risk of injury~\cite{OesthLarssonJakobsson22}.
Muscle is an incompressible tissue~\cite{Swammerdam1758} that is driven by a vast array of microscopic contractile fibers. 
Models of muscle often ignore the volume of muscle, and instead is treated as an array of massless contractile cables.  
The tension that muscle can generate is challenging to model because the stress, stiffness~\cite{SugiTsuchiya88}, and damping~\cite{KirschEtAl94} all vary with chemical activation, strain~\cite{GordonEtAl66}, strain-rate~\cite{Hill38}, and state history~\cite{HerzogLeonard00}.

\begin{figure}[tbp]
	\centerline{\includegraphics[width=17.5cm]{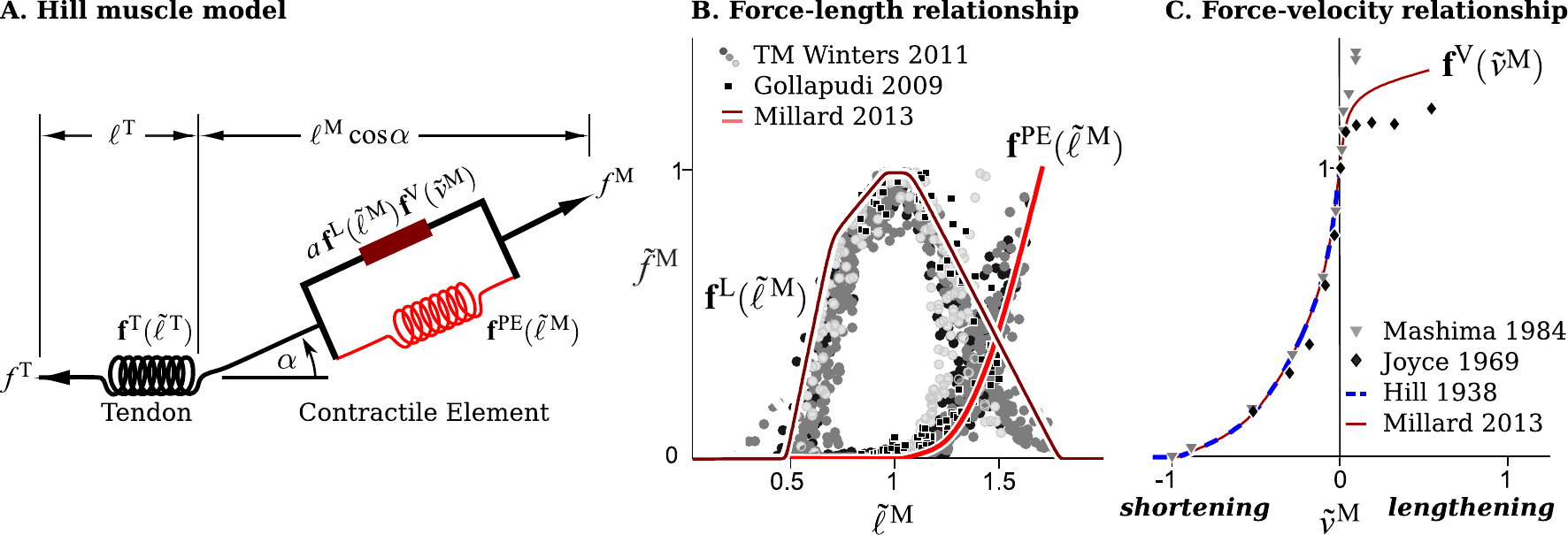}}
	\caption{\enspace The most commonly used muscle model is a phenomenological Hill-type muscle model (A). A Hill-type muscle has two main parts: an elastic tendon, and a contractile element that acts at a pennation angle $\alpha$ to the tendon. When normalized force $\fmN$ is plotted against the normalized length of the muscle $\lmN$ the nonlinear and scale-invariant force-length relationship arises~\cite{WintersEtAl11, GollapudiLin09}. Similarly, when normalized muscle force is plotted against normalized contraction velocity $\vmN$ (C) Hill's hyperbolic force-velocity relationship appears~\cite{Hill38}. During lengthening, the pattern of force development is not as consistent, but usually involves enhanced forces~\cite{Joyce69,Mashima84}.  Splines are often used to represent the nonlinear force-length $\flN$ and force-velocity $\fvN$ relationships~\cite{MillardEtAl2013}. \label{Fig:Muscle}}
\end{figure}

The most commonly used muscle model is the phenomenological Hill-type model~\cite{Hill38,Wilkie56,RitchieWilkie58} because it strikes a balance between accuracy and computational cost.
Modern Hill-type models~\cite{Zajac89,MillardEtAl2013} are formulated in a normalized space\footnote{Note that normalization is indicated with a tilde, forces are normalized by $\fopt$, CE lengths by $\lopt$, CE velocities by $\vmax$, and tendon lengths by $\ltslk$.} which allows a Hill model to be fit to any musculotendon actuator by knowing only a few architectural parameters of the biological muscle: the maximum isometric force $\fopt$, the length $\lopt$ at which the CE develops $\fopt$, the pennation angle $\pen$ between the CE and the tendon at $\lopt$, the maximum shortening velocity $\vmax$ of the CE, and the slack length $\ltslk$ of the tendon. 
Embedded in the formulation of the Hill model are three experimental characteristics: Hill's~\cite{Hill38} force-velocity relation $\fvN$, Gordon et al.'s active-force-length $\flN$ relation~\cite{GordonEtAl66}, the passive force-length relation $\fpeN$ of the CE~\cite{GollapudiLin09,WintersEtAl11}, and the nonlinear force-length relation $\ftN$ of the tendon~\cite{MagnussonEtAl01, MaganarisPaul02}.
The normalized tension developed by a Hill-type CE is the sum of the active and passive forces 
\begin{equation}
	\fmN =  \act \, \flN \, \fvN + \fpeN 
\end{equation}
developed by the CE, where $\act$ is the chemical activation of the muscle.
The activation of the muscle is driven by the input signal $u$ of the nervous system
\begin{equation}
\dact = \left\{\begin{array}{lr}
(u-\act)/(\tcA\,(0.5+1.5\act)) & \text{ if } u > a,\\
(u-\act)/(\tcD/(0.5+1.5\act)) & \text{otherwise}
\end{array}\right.
\end{equation}
where the time constants $\tcA$ and $\tcD$ that are largely determined by rates at which calcium ions can be pumped into and out of sarcomeres, respectively.
The tendon is often modeled as a nonlinear massless elastic~\cite{Zajac89, MillardEtAl2013} that is assumed to be in a perfect force equilibrium
\begin{equation}
	\fmN \cos \pen - \ftN = 0 \label{eqn:HillEqEqn}
\end{equation}
with the CE.
This force equilibrium equation can be solved for $\vmN$ by inverting $\fvN$, or by regularizing this equation with a numerical damping term and applying a root solving method to solve for a $\vmN$ that satisfies Eqn.~\ref{eqn:HillEqEqn}~\cite{MillardEtAl2013}.
The final result is a first order ODE that has two states $(\act,\lmN)$ and can replicate Hill's~\cite{Hill38} and Gordon et al.'s~\cite{GordonEtAl66} experiments by construction.
Often large muscles are represented in MBS and FE models by an array of Hill-type contractile elements that attach to, and wrap around, bony geometry~\cite{ScholzEtAl16}.

Once a model is developed, there are two different families of methods that can be used to simulate multibody~\eqref{Eq:Mubo_equation} and FE systems~\eqref{eq:Nonlinear_Fe_Equation}: inverse methods and forwards methods.
Although both inverse and forwards methods can be applied to multibody~\eqref{Eq:Mubo_equation} and the FE systems~\eqref{eq:Nonlinear_Fe_Equation}, FE systems are almost exclusively simulated using forward dynamics.
There are two inverse methods, inverse kinematics and inverse dynamics, that are often applied to MBS systems.
The inverse kinematic method is frequently used with HBM models when analyzing experimental data. 
Unfortunately it is rarely possible to measure joint angles in a human participant, instead the positions $\mkr$ of motion capture markers placed on the skin are measured using camera systems.
The generalized positions $\yb$ of the model that best fit observed motion are solved for by minimizing the sum of squared errors between the $\mathrm{P}$ recorded marker positions $\mkr$ and the model's virtual marker positions $\vmkr(\yb)$
\begin{equation}
\min \sum_i^\mathrm{P} \left( \mkr_i - \vmkr_i(\yb) \right)^{\mathrm{T}}\cdot\left( \mkr_i - \vmkr_i(\yb) \right)
\end{equation}
using the Levenberg-Marquardt algorithm~\cite{Levenberg44,Marquardt63}.
By repeating this process at every point in time, the full trajectory of a model $\yb(t)$ can be constructed and numerically differentiated to yield generalized velocities $\ybp(t)$ and accelerations $\ybpp(t)$.
Inverse dynamics analysis is an additional method that uses $(\yb,\ybp,\ybpp)$ and Eqn.~\ref{Eq:Mubo_equation} to solve for the external forces ${\boldsymbol{k}}_{\mathrm{e}}$ that must have been applied to the model during the experiment.
In the context of driving, inverse dynamics is often used to solve for the net muscle forces developed in the leg during a breaking maneuver.
In addition, inverse dynamics can also be used to evaluate Eqn.~\ref{Eq:Mubo_equation} as a constraint using synthetically generated $(\yb,\ybp,\ybpp)$ as is done in the direct collocation method~\cite{VonStryk93}.
Since inverse dynamics does not require any computationally expensive numerical operations, such as matrix inversion or numerical integration, it is often used in optimization, parameter identification, and in the analysis of pre-recorded human motion.

In contrast to inverse methods, forward dynamics can be used to predict the trajectory that a model will take over time given only its initial state and a function to evaluate the forces applied to the system.
Using the FE system as an example, the forward dynamic method begins with the initial positions $\qb_0$ and velocities $\qbp_0$ 
\begin{align} 
	\begin{split}
		\Mb \cdot \qbpp+ \fb_\mathrm{int}(\qb,\qbp)&=\fb_\mathrm{ext}\\
		\text{with}\hspace{1cm}t &\in [ t_0, t_{\mathrm{end}}],\\
		\qb(t_0) &=\qb_0,\\
		\qbp(t_0) &=\qbp_0,\\
		\qb &= \hat{\qb} \text{ on } \partial \mathrm{B}_q,\\
		\bm{\sigma} &= \hat{\tb} \text{ on } \partial \mathrm{B}_{\sigma}
	\end{split}
\end{align}
and a numerical integration scheme that steps forwards in time from $[ t_0, t_{\mathrm{end}}]$ using the value of $\qbpp$ from the equations of motion. 
The boundary conditions\ $\qb=\hat{\qb} \text{ on }\partial \mathrm{B}_q$, and the initial stress $\bm{\sigma}=\hat{\tb} \text{ on } \partial \mathrm{B}_{\sigma}$, as well as the initial positions are often crucial for valid and realistic simulation results. 
The choice of suitable integration method for numerical simulation depends on the type of FE equation, smoothness of the data, and the response of interest. For highly-dynamic impact situations, explicit integration schemes are preferred and state-of the art, since the penalty-based contact algorithms introduce noise~\cite{BelytschkoLiuMoran00}.

\begin{figure}[t]
	\centerline{\includegraphics[scale=1]{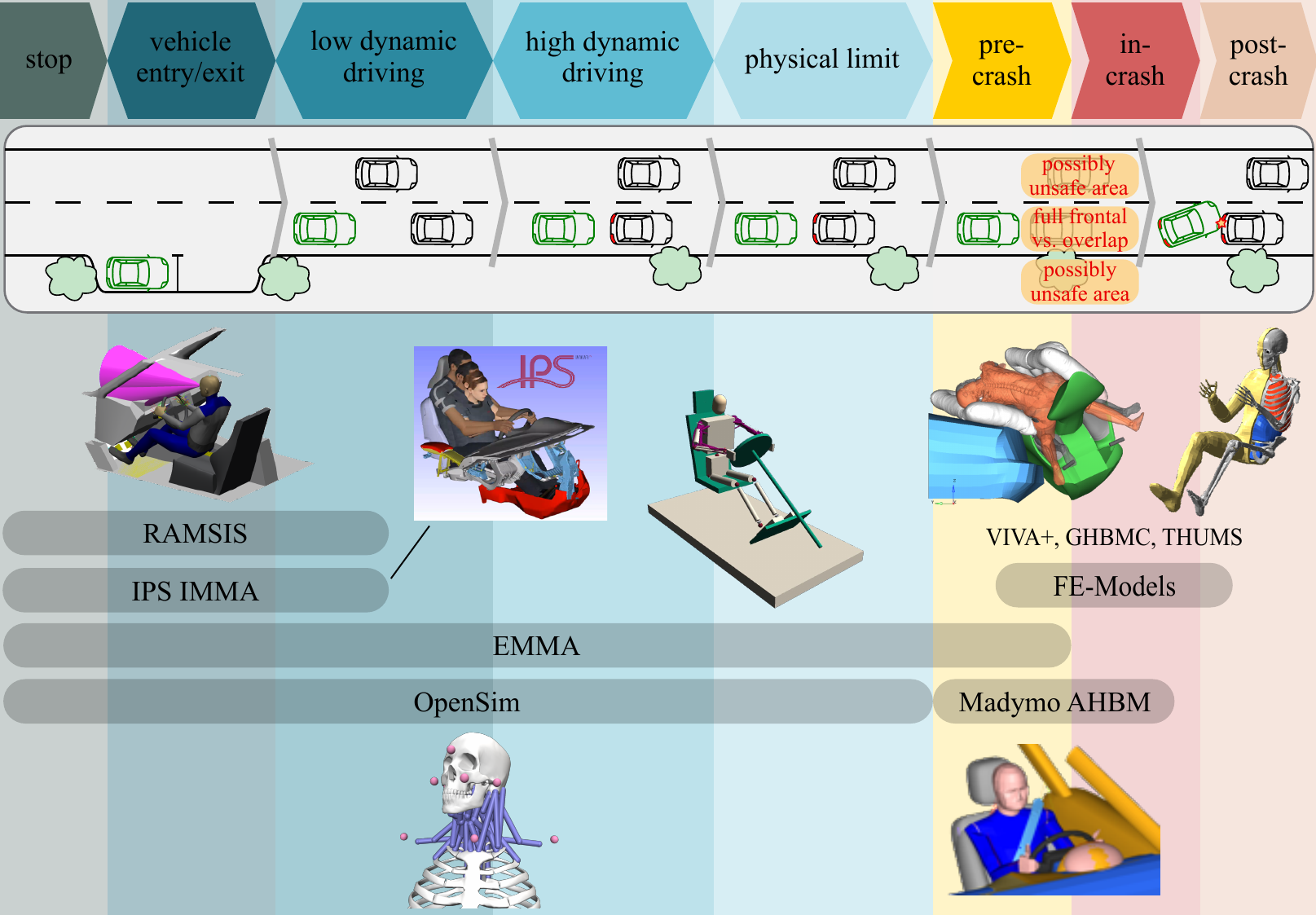}}
	\caption{\enspace We have defined a series of driving phases to guide us to a model that is just detailed enough to answer typical questions associated with each phase. While the stop and entry/exit phases are self-evident, low-dynamic driving is defined by low accelerations: the vehicle might be traveling slowly or at highway speeds, but the acceleration is negligible. During these phases, kinematic or multibody models can be used. The vehicle accelerates during high-dynamic driving but remains within its physical limits. The vehicle is in the physical-limit phase once it can no longer follow the desired trajectory. A multibody model should be used to incorporate the influence of acceleration on the body during these phases. The in-crash phase begins once a collision is imminent, transitions to the in-crash phase once contact is made, and ends with the post-crash phase once the colliding bodies have come to rest. Finite element models should be used to simulate the in-crash phase because the accelerations are so high that the body's segments cannot be approximated as rigid. The abbreviations of the HBMs appear in Sec.~\ref{sec:human_body_models_vehicle_safety}.\label{fig_DrivingSituations}}
\end{figure}

We have found it helpful to define a series of driving phases (upper part Fig.~\ref{fig_DrivingSituations}) which apply increasingly large loads to the body: vehicle entry/exit, low dynamic driving, high dynamic driving, physical limit, pre-crash, in-crash, and post-crash. 
Vehicle entry/exit occurs when passengers enter or exit a static vehicle.
During low-dynamic driving, the vehicle moves, but the acceleration is so small that it can be ignored.
Kinematic or multibody models might be used within these phases, as ergonomics is often the primary concern.
Beyond the low-dynamic-driving phase, a multibody model must be used to capture the body's response to acceleration.

The high-dynamic driving phase occurs when the car accelerates within its physical limits.
The physical-limit phase is reached when a desired trajectory can no longer be followed because of a physical limitation, such as a wheel slipping or a loss of wheel-ground contact.
During these phases, a multibody model is often sufficient because the accelerations are high enough to move the passenger's body but are not significant enough to prevent voluntary motion or cause injury.

Once a collision is imminent, the pre-crash phase has begun and does not end until contact is made between the vehicle and the colliding body. 
The in-crash phase begins with the colliding body contacts the vehicle and ends when the vehicle comes to rest, after which the post-crash phase begins. 
The accelerations and forces applied to the body during the in-crash phase can be high enough that a finite-element model is needed to simulate the stresses, strains, and injuries to the body. 
Note that we have modified Haddon's ~\cite{Haddon70} original definitions of the crash phases to define these phases based on the kinematics of colliding vehicles and contact events to make these definitions more applicable to simulation.

The final three phases are all associated with a crash and can cause severe injury: the pre-crash phase, which consists of the seconds before impact; the in-crash phase, which often lasts only 100-200 ms; and the post-crash phase, which begins after the vehicle has come to rest.
The accelerations applied to the body during the in-crash phase can be quite high: during a frontal crash against a non-deformable barrier, the occupant of the car can be decelerated from $v^{\mathrm{occ}}_{\mathrm{in-crash}}=43\;\mathrm{km}/\mathrm{h}$ to $v^{\mathrm{occ}}_{\mathrm{post-crash}}=0 \;\mathrm{km}/\mathrm{h}$ within~105\;ms~\cite{ProchowskiEtAl11}, which can cause severe injury and death.
While the severity of injury can be evaluated from the stresses and strains of the body's tissues, medical and economic institutions classify injury differently: the medical \emph{Abbreviated Injury Scale} (AIS)~\cite{AAAM08} rates injuries of the whole body from 'non-injured' (0) to 'untreatable' (6); in contrast, economists value injury in terms of the years of life lost. 
Accordingly, national traffic safety organizations and insurance companies carefully monitor mortality rates and the subsequent costs due to road traffic.

The paper continues with an overview of available HBMs and a series of case studies. 
Our most comprehensive example in Sec.~\ref{sec_application} presents an analysis of a novel motorcycle safety system using: 
(i)~in-crash FE simulations and injury evaluation of riders on a motorcycle equipped with a novel safety system;
(ii)~validation experiments using a Driver-in-the Loop system to access the accuracy of a Madymo multibody HBMs pre-crash behavior;
(iii)~a prediction of human behavior during a take-over scenario using the muscle-actuated EMMA model~\cite{RollerEtAl20}.
The paper concludes with current research trends in vehicle safety simulation.

\section{Human Body Models in Vehicle Safety}\label{sec:human_body_models_vehicle_safety}

A variety of HBMs have been developed to strike a balance between model detail and simulation runtime for the wide range of automotive design concerns.
Tab.~\ref{tab_overviewTable} shows how the level of detail varies widely across the HBMs we consider, from models that are kinematic and driven by tabulated experimental data to models that include accurate representations of detailed anatomical structures of the body.
The level of detail in the model directly translates into a greater memory footprint, larger data structures, and longer mathematical operations during simulation: detail is resource intensive.
As such, it is most common to use kinematic and multibody models to solve optimization problems that may require hundreds to thousands of simulations of a specific scenario.
In contrast, when high levels of detail are needed an FE model is typically used and the crash scenario is often only simulated a few times.
Model choice is often guided by the level of detail required for the study, and limited by the computational costs of running the simulations.

\subsection{Kinematic Human Body Models}
At the lowest end of modeling complexity, kinematic models are mainly used to evaluate ergonomics.
These models only consider the kinematic motion variables and introduce joint angle limits to evaluate reachability.
Examples for this class of HBMs are RAMSIS~\cite{BubbEtAl06} and Jack~\cite{PhillipsBadler88} which are both connected to anthropometric databases to evaluate the ergonomics of a specific age, sex, or population.
It is further possible to visualize reach envelopes and to evaluate the visibility of instruments, mirrors, and the environment~\cite{BubbEtAl06, PhillipsBadler88}.
In addition, RAMSIS can generate realistic postures and movements with a probability-based model that relies on pre-recorded posture studies and can evaluate the discomfort of these postures~\cite{BubbEtAl06}.
RAMSIS is widely used in the automotive industry to evaluate the interaction between the human and product during the early design phase~\cite{Wirsching19}.

\subsection{Multibody Human Body Models}

Multibody HBMs offer a balance between the level of detail of the model and the computational cost of simulation.
This balance makes it possible to use multibody models in many different phases of driving from vehicle entry/exit, ergonomic simulations, low-dynamic driving, high-dynamic driving simulations and pre-crash simulations.
In addition, the level of detail in a multibody model can vary dramatically, from HBMs that are driven by joint torques and treat contacts as rigid constraints~\cite{BjoerkenstamEtAl18}, to HBMs that are driven by hundreds of lumped-parameter muscles and approximate contact forces using compliant contact volumes~\cite{SiemensMadymoModelManual20}. 

\paragraph{OpenSim}
OpenSim is an open-source software system that allows users to simulate a wide variety of movements of musculoskeletal models~\cite{DelpEtAl07} that can vary from simple joint torque driven models to models that contain hundreds of nonlinear muscle actuators.
It is possible to interface with the OpenSim software via Python, Matlab and C++, among others.
OpenSim includes functionality for forward/inverse kinematics, forward/inverse dynamics, and optimal control~\cite{DelpEtAl07, SethEtAl18} making it a tool that can be used to analyze complex experimental data, or synthesize physics-informed human movement predictions.
Most often OpenSim is used to compute physics-based estimates of the muscle forces and joint contact loads from an experimental recording of a human participant.
In~\cite{IndraniHakanssonLankarani19} OpenSim is used to evaluate the car-pedestrian impact.

\paragraph{AnyBody}
AnyBody~\cite{DamsgaardEtAl06} provides similar dynamic and kinematic functionality compared to OpenSim but focuses exclusively on solving an inverse problem: estimating the muscle forces and joint loads of a participant from motion capture data of an experiment. 
As long as the model and input data are identical, both OpenSim and AnyBody produce the same results~\cite{KimEtAl18}.
However, because AnyBody has focused on solving inverse problems, they have been able to make incredibly detailed musculoskeletal models that include hundreds of muscles and detailed joint models.
This additional level of detail does result in differences when simulations of AnyBody models and OpenSim models are compared~\cite{KimEtAl18, SandholmEtAl11, WagnerEtAl13}.

\begin{figure}[tbp]
	\centerline{\includegraphics[width=1.0\textwidth]{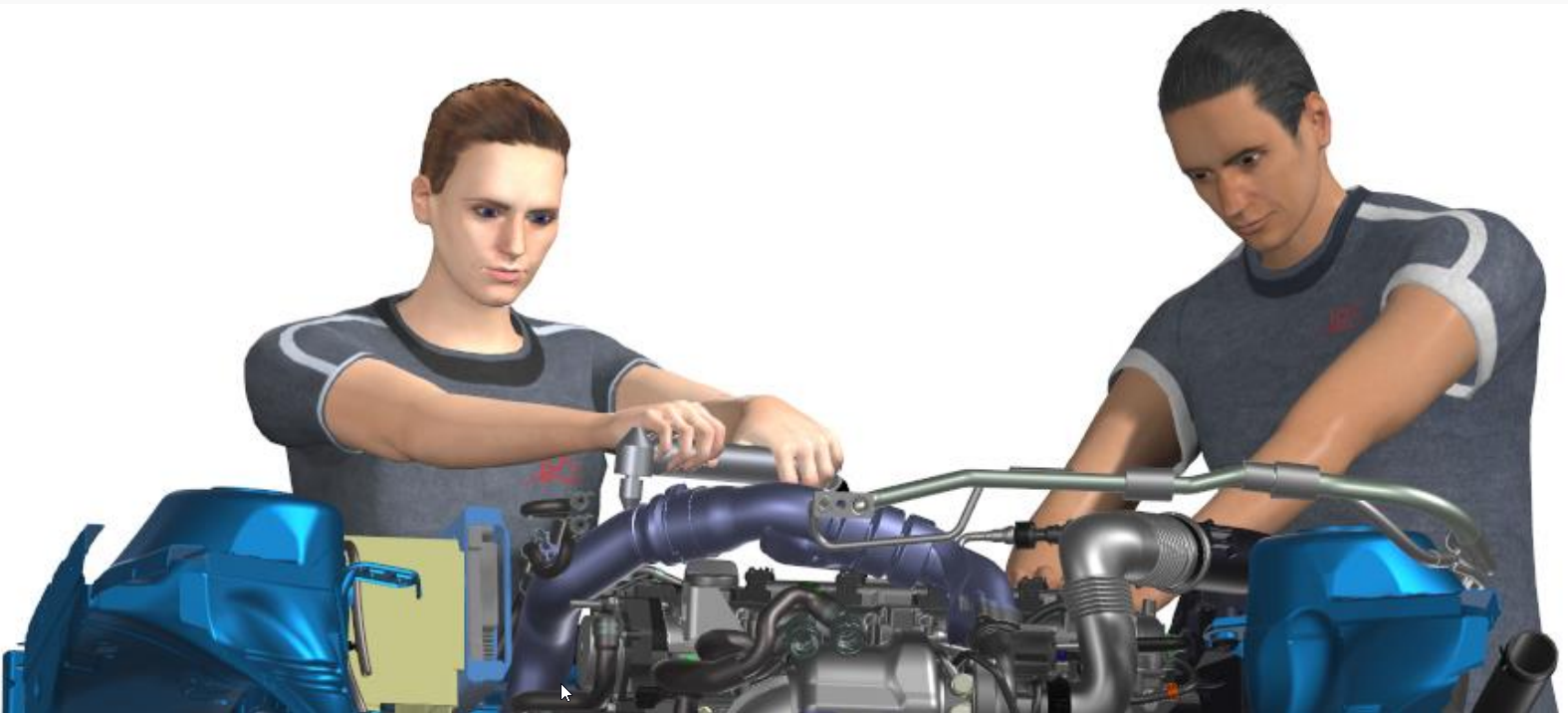}}	
	\caption{\enspace Simulation of the assembly of an engine of two workers in IPS IMMA .\label{fig_IMMA_2worker}}
\end{figure}

\begin{figure}[tbp]
\begin{minipage}[t]{0.48\textwidth}
	\centerline{\includegraphics[width=1.0\textwidth]{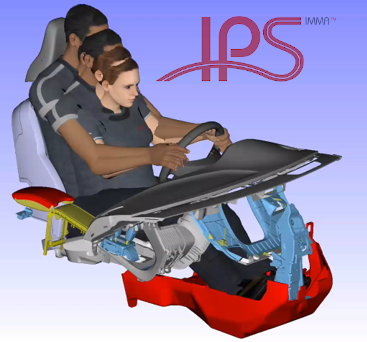}}	
	\caption{\enspace Simulation of a family of DHM sitting in a car holding the steering wheel.\label{fig_IPS_VD}}
\end{minipage}
\hfill
\begin{minipage}[t]{0.48\textwidth}		
	\centerline{\includegraphics[width=1.0\textwidth]{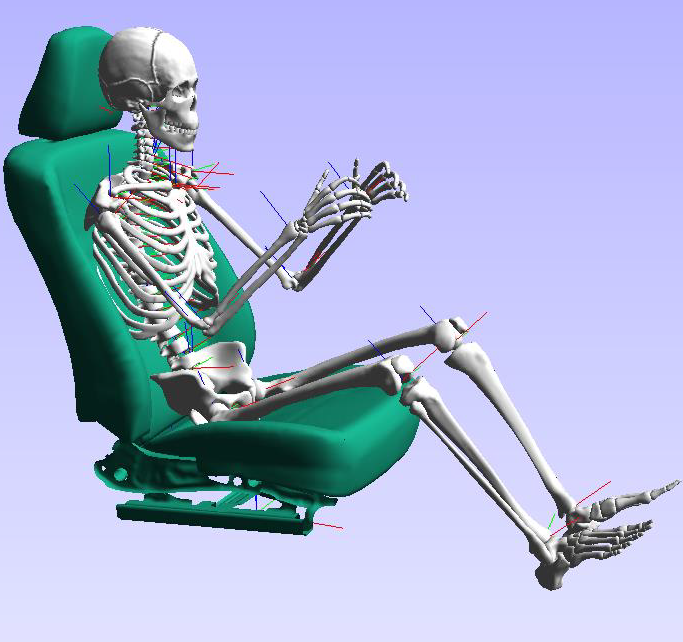}}	
	\caption{\enspace The skeleton of the EMMA4Drive manikin based on the THUMS Skeleton.\label{fig_EMMA4DRIVE}}
\end{minipage}
\end{figure}

\paragraph{IPS IMMA}
The Intelligently Moving Manikin (IMMA)~\cite{HansonEtAl19} is an MBS focused on evaluating human ergonomics~\cite{HogbergEtAl18} that has been applied to study human tasks in a production environment, and also to study the ergonomics of seated reaching tasks in a vehicle.
IPS~IMMA is a biomechanical model that consists of 82 bone segments connected by joints, giving a total of 162 DOFs. 
The model creates a human-like posture due based on these constraints, by solving an optimization problem.
This includes finding a stable equilibrium (balance),  which minimizes a proprietary objective function that contains terms related to ergonomics, joint torque, and posture.
The posture of the manikin is controlled by the user defines by defining where to grasp an object and/or where the manikin should look, without the need to touch any joint angles.
IPS~IMMA creates a motion using a sequence of optimized static postures (quasi-static simulation).
These movements are specified by defining a series of boundary poses in a sequence editor.
Using the sequence of boundary conditions as a guide, the IPS~IMMA system computes a series of postures interpolating between the sequence of poses that predicts quasi-static human movements.
IPS~IMMA also includes a reachability analysis and mountability feature, which ensures collision-free assembly movements.
Depending on the task, the movement strategy is adjusted automatically.
This feature also allows for the examination of the mountability of objects by humans.
As manikins are controlled by defining constraints, which remain the same for manikins with another anthropometry, a defined scenario can be recalculated for a whole manikin family without any further control effort.

In recent days IPS~IMMA was extended to simulate manikins sitting in a car or a truck, see Fig.~\ref{fig_IPS_VD}.
Herein, different objective functions are used to create the static posture tailored to these use cases.
Therefore, IPS~IMMA can also be used for static reachability analysis or visibility studies like RAMSIS.
For dynamic use cases, the EMMA model has been developed in recent years.

\subparagraph{EMMA}

EMMA (Ergo-dynamic Moving Manikin) is a multibody model that has been used to synthesize human motions when coupled with an optimal control solver and a cost function.
The movements of EMMA are induced by Hill-type muscles~\cite{Hill38, MaasLeyendecker13} which are represented as cables, or by joint actuators which act directly on the joints of the model.
EMMA's synthesized movements are driven to minimize cost functions that consist of several terms including the sum of squared control signals, the kinetic energy of the system, and the execution time of the movement.
Although optimal control problems applied to musculoskeletal systems can be difficult to solve, long scenarios have been simulated using a specialized optimal control problem formulation, \textit{discrete mechanics and optimal control for constrained systems} (DMOCC)~\cite{LeyendeckerEtAl10}, and an interior point optimization method~\cite{BjoerkenstamEtAl18, RollerEtAl17}.

EMMA uses a special case to solve optimal control problems that include contact: if the contact sequence is known ahead of time the optimal control problem can be broken up into a series of discrete phases.
Each discrete phase of the optimal control problem is defined by the set of bodies that are in contact, since the equations of motion need to be modified to include the contact.
As is commonly done when simulating walking~\cite{Koch12} and running~\cite{Schultz09}, contact is treated as a kinematic constraint with additional phase-specific constraints added to the optimal control problem to ensure that the contact does not exceed friction limits.
Transitions between phases are made assuming that contact is plastic, and momentum conserving. 
As such, every phase transition is accomplished by applying the principle of impulse momentum to evaluate the generalized velocity vector of the system after each contact. 
Finally, the duration of each contact is left as a free variable that the solver modifies to improve the cost of the solution.
Since the duration of a phase can even go to zero, effectively eliminating the phase, this approach can be applied to a wide variety of problems.

Typically, a multi-phase optimal control problem is limited by the duration and complexity of the task being simulated.
However, because EMMA makes use of a discrete mechanics formulation DMOCC~\cite{LeyendeckerEtAl10}, the physical properties (impulse, angular momentum) of the system are conserved even for large discretization steps, which allows EMMA to simulate tasks that have a long duration.
Thus far EMMA has been used to simulate a person lifting a box~\cite{RollerEtAl17}, pulling an emergency brake in a car~\cite{RollerEtAl20}, and getting into a truck~\cite{BjoerkenstamEtAl20}.
EMMA is currently being used to simulate the voluntary motions that occur just prior to the start of a crash, bridging the gap between MBS and FE models, see Fig.~\ref{fig_EMMA4DRIVE}. 

\paragraph{Simcenter Madymo}

Simcenter Madymo is a combined MBS and FE simulation system that features full-scale omnidirectional multibody models~\cite{HappeeEtAl98}.
There are a number of advanced elements that distinguish Simcenter Madymo such as elastic ellipsoid contact models, detailed facet models represented with meshes, flexible components such as ribs, and FE elements to represent important deformable elements~\cite{HappeeEtAl03}. 
The focus of Simcenter Madymo is the simulation of active muscle-induced human motion during the pre-crash phase, and as such, the models are referred to as Active Human Body Models (AHBMs).

\begin{figure}[tbp]
	\begin{minipage}[t]{0.48\textwidth}
		\centerline{\includegraphics[width=1.0\textwidth]{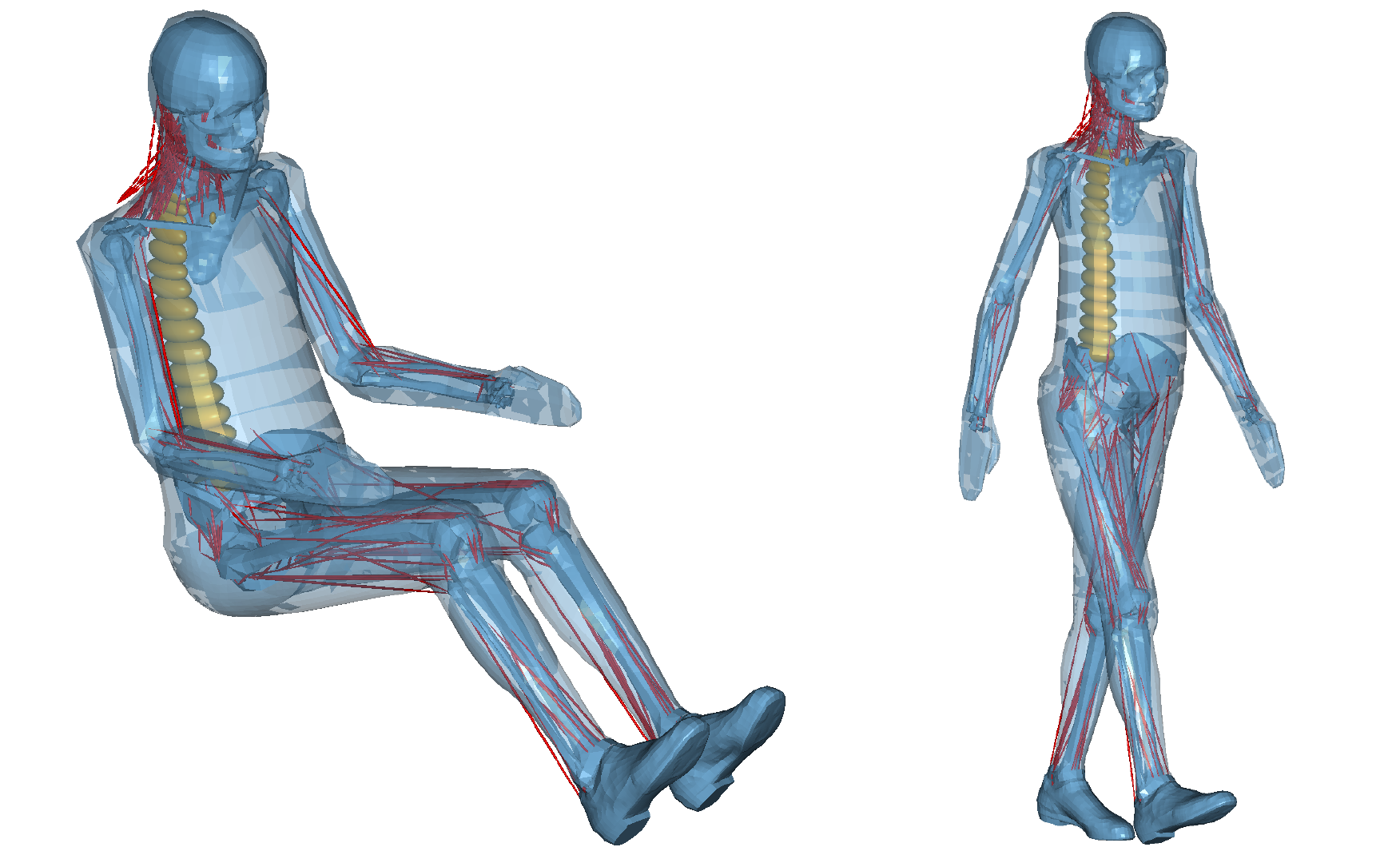}}	
		\caption{\enspace Active Human Body Model in the occupant (left) and in the pedestrian/standing version (right). The red lines correspond to the muscle elements. The spine is modeled as flexible elements~\cite{SiemensMadymoModelManual20}.\label{fig_AHBM_Pedestrian_Occupant}}
	\end{minipage}
	\hfill
	\begin{minipage}[t]{0.48\textwidth}		
		\centerline{\includegraphics[scale=1.0]{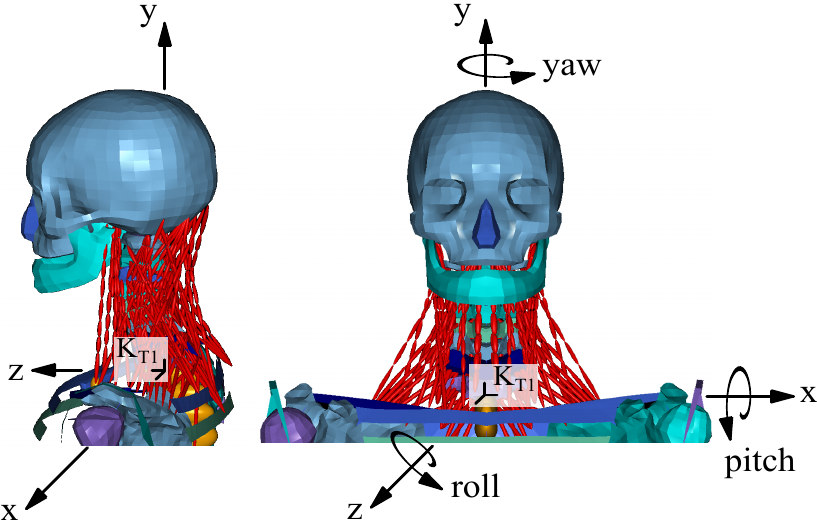}}
		\caption{\enspace Frontal and side view of the AHBM head and neck model showing 32 (2x16) different muscles modeled with 136 muscle elements. The rotation angles are also depicted.\label{fig_HN_Model_Overall}}
	\end{minipage}
\end{figure}

The AHBM is a detailed model of a 50\% male from the RAMSIS database (a $75.3\;\meh{kg}$ and $1.76\;\meh{m}$) that features a seated (Fig.~\ref{fig_AHBM_Pedestrian_Occupant} left) and a standing variant (Fig.~\ref{fig_AHBM_Pedestrian_Occupant} right) with the body being supported by a skeleton that has 182 rigid bodies, 8 flexible bodies (abdomen and thorax) and compliant joints~\cite{SiemensMadymoModelManual20}.
The AHBM is not fixed, but supports modification of the joints and surface meshes of the model.
The outer surface of the model uses FE shells to accurately model contact forces with the seat, interior, and safety systems. 
Effort has been put into defining the contact stiffness of different parts of the vehicle and the AHBM to improve the realism of contact simulation.

\begin{figure}[tbh]
	\centerline{\includegraphics[scale=0.9]{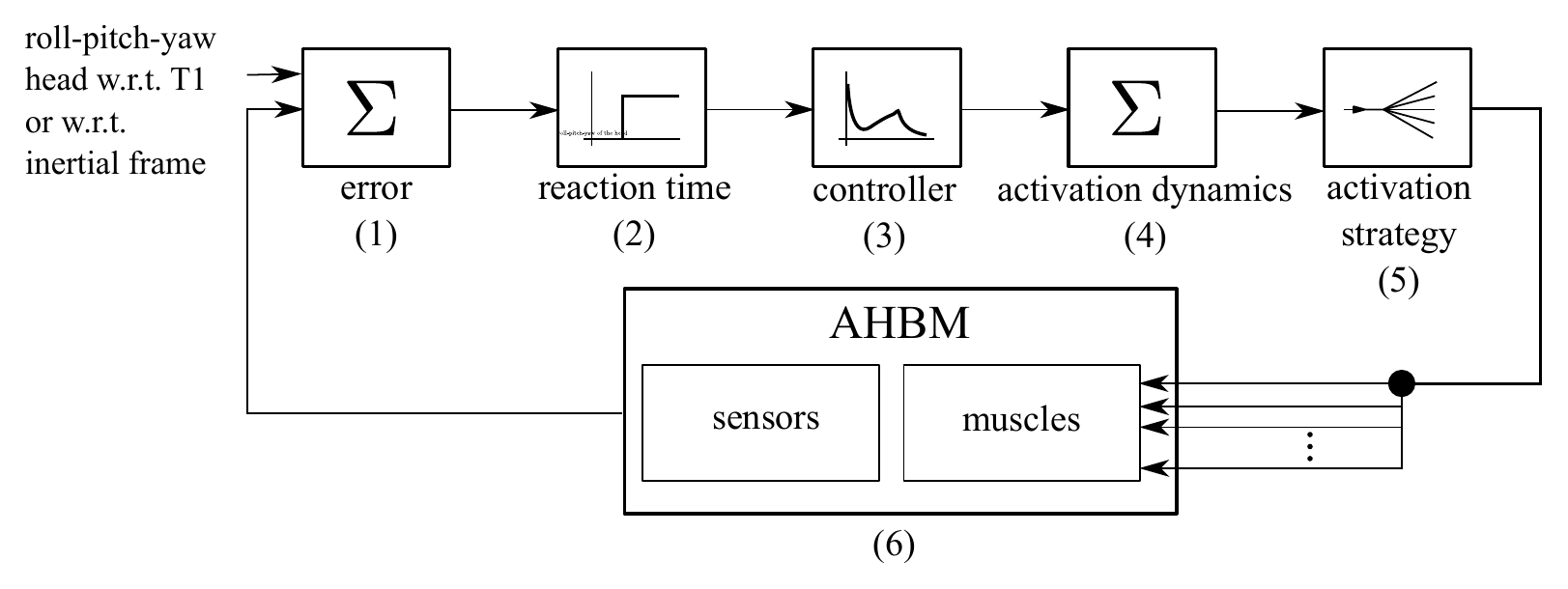}}
	\caption{\enspace Control approach of the Madymo AHBM. 
		\label{fig_madymo_ahbm_regelkreis}}
\end{figure}

Two unique features of Madymo's AHBM are its posture stabilization control system, and its anatomically detailed active musculature of the extremities (Fig.~\ref{fig_AHBM_Pedestrian_Occupant}) and neck (Fig.~\ref{fig_HN_Model_Overall}).
The posture stabilization system consists of a set of reflex controllers that have been carefully designed to mimic the reflex responses of the human body across a range of bandwidths.
These reflex controllers are used to activate the muscles  (Fig.~\ref{fig_madymo_ahbm_regelkreis}) of the head and neck, arms, and legs, while the spine is driven by torque actuators.
While all of the controllers are important, the controller for the head and neck is particularly important because the cervical spine is frequently injured~\cite{Pink16} during car accidents.

Both the musculoskeletal model (Fig.~\ref{fig_HN_Model_Overall}) and control of the head and neck is complex (Fig.~\ref{fig_madymo_ahbm_regelkreis}).
While the pose of the model's head is described by only three coordinates (pitch, roll, and yaw), the neck and head are driven by 136 muscle elements belonging to 32 different muscle groups that represent the most relevant muscles of the head and neck.
Reflex and stabilization studies of the head and neck~\cite{ForbesEtAl13, HappeeEtAl17} have shown that the neck's muscular recruitment is driven by input to the somatosensory system, and that the response is frequency-dependent.
For static poses and low frequency movements, the vestibular system is used to achieve an upright position of the head relative to the inertial system. 
However, in more dynamic situations with higher frequency content, the system tries to control the relative pose between the head and torso (head-in-space stabilization).

Accordingly, the AHBM uses this postural control concept and tries to assume a specific pose as quickly as possible under the influence of disturbing forces acting on the body.
Reflexes are incorporated into the head and neck controller (HNC)~\cite{FragaEtAl09, NemirovskiVanRooij10} using a model of the somatosensory system to trigger reflexes, which are subject to a conduction delay, and are modified by a frequency-dependent gain prior to activating the muscles (Fig.~\ref{fig_madymo_ahbm_regelkreis}).

The input to the HNC system (Fig.~\ref{fig_HN_Model_Overall}) is the desired orientation trajectory of the head which is defined as the orientation of the head in the inertial system (using roll-pitch-yaw), and the orientation of the head with respect to the T1 vertebrae (also using roll-pitch-yaw) which are delayed prior to reaching the muscles.
The extra coordinate system allows the system to smoothly switch between using an inertial system orientation reference and a head-in-space orientation reference.
The error signal between the desired and actual head orientation is delayed before reaching the muscles.
As the signal path is shorter to the neck than to the extremities, the transmission delay for activating the muscles of the neck is shorter than for the arms or legs.

Different muscle groups are stimulated to counteract the deflection of the head and neck, and to stabilize the movement of the head using co-contraction.
With its 32 muscles, discretized to 136 muscle elements, the AHBM has a highly over-actuated head and neck.
A simplified linear and decoupled control system is used(Fig.~\ref{fig_madymo_ahbm_regelkreis})
The control system consists of three time-delayed PID controllers and a co-contraction factor $\textit{CCR}_{\textrm{Neck}}$.
A higher co-contraction level makes the HNC less sensitive to external disruptive forces.
A spatial activation pattern derived from experimental data~\cite{OesthEtAl17, LarssonEtAl19} is used to activate the over-actuated neck model.
Madymo's control strategy is not hard-coded but can be customized by using the Matlab Madymo coupling strategy.
For example, Kempter et al.~\cite{KempterBechlerFehr20} implemented a completely different control strategy in which every individual muscle of the neck was controlled by its own stretch-reflex controller.
By using the Matlab-Madymo coupling strategy, it is possible to update the HNC as more about this complex biological system is learned.
Many properties of the AHBM model have been experimentally validated, and calibrated.
The passive properties of the AHBM have been validated using PMHS tests~\cite{SiemensMadymoModelManual20}.
Many active properties of the model have been indirectly validated by comparing the response of the entire head and neck model to in-vivo sled tests~\cite{ChoiEtAl05, EjimaEtAl07, KirschbichlerEtAl11}, Driver-in-the-Loop (DiL) tests~\cite{FehrEtAl21b, KempterEtAl18, KempterEtAl23}, and tests in real vehicles~\cite{HuberEtAl15, OesthEtAl12}.
An overview of available tests is provided by the OSCCAR consortium~\cite{OSCCAR22}.

\subsection{Passive Finite Element Models}

The FE method has been used extensively~\cite{YangEtAl06} to create HBMs that accurately capture the geometry and injury mechanics of the musculoskeletal system, the nervous system, and internal organs of the body.
Passive FE models were originally developed to mimic the response of a PMHS to occupant restraint systems.
The word passive is used to denote the fact that these models do not include active muscles, nor any active actuation.
Even though muscles remain inactive in a passive FE model, these models are sophisticated and can predict injury to bones, connective tissue, and internal organs.
A broad selection of element formulations have been developed to strike a balance between accuracy and runtime during simulation.

Element formulations range in terms of geometry (hexahedral and tetrahedral are typical), computational properties (DOF, numerical features), the physical interactions (large deformations, contact), and injury mechanisms (fracture, rupture, tearing) that are supported.
The geometry of the structures in an HBM is often based on CT or MRI images and contains many irregular shapes that are difficult to mesh.
Even though irregular shapes can be easily meshed using tetrahedral elements, the preferred element for HBMs is the hexahedral element due to its superior numerical robustness during simulation.
The formulation of the element, however, is not the only modeling component that affects simulation time.

A variety of joint models are used in HBM models to capture the desired level of detail.
Idealized joints, such as spherical and revolute joints, are often used to reduce computation when injury to the joint is not the focus of the simulation. %
In contrast, anatomically accurate joint models are used when injury to the cartilage, bursa, and ligaments of a joint are of interest. %
A wide variety of models have been developed to balance detail against simulation time for the many crash scenarios that must be considered during the design phase.

\begin{figure}[tbp]
	\centerline{\includegraphics[width=17.3cm]{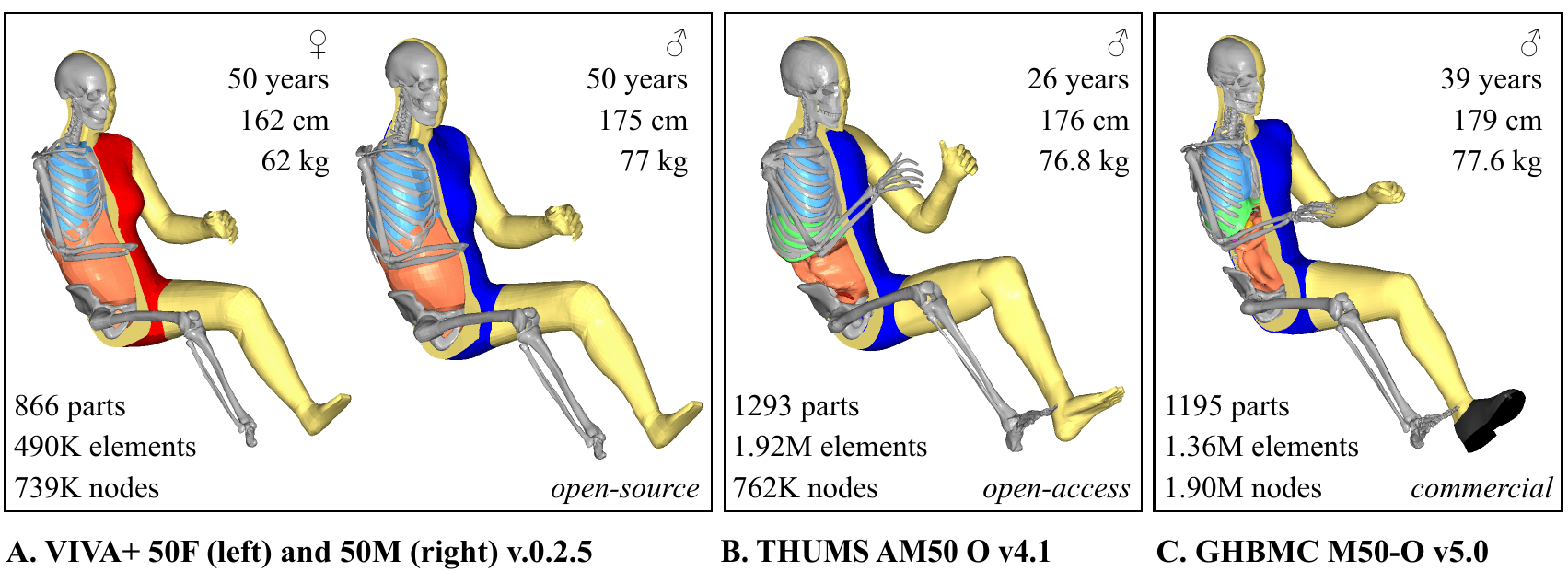}}
	\caption{\enspace Selection of female and male HBMs~(partially blanked) in their default positions~(\cite{MaierEtAl22}). \label{Fig:SafeMoto_HBMDefault}}
\end{figure}

Current state-of-the-art HBMs in vehicle safety (Fig.~\ref{Fig:SafeMoto_HBMDefault}) include Toyota's Total HUman Model for Safety (THUMS)~\cite{KatoEtAl18}, the VIVA+~\cite{JohnEtAl22} models, and the GHBMC~\cite{SchwartzEtAl15, GHBMC22, CorreiaMclachlinCronin21} models. 
Six THUMS model versions have been developed and evaluated~\cite{Toyota21B} since 2000.
VIVA+ is a family of open-source models that have been developed as a part of the EU-VIRTUAL consortium~\cite{LindnerEtAl20}. 
These models are derived from the VIVA+ OpenHBM --- a 50$^{th}$ percentile female model from Chalmers University, Gothenburg, Sweden --- that was developed to study Whiplash Associated Disorders (WAD) of women during rear-end collisions~\cite{OesthEtAl17a}.
The Global Human Body Models Consortium (GHBMC) is composed of North American car manufacturers and research institutions. 
Using MRI and CT data of a 26 year old male, the GHBMC developed its first $50^{th}$ percentile male HBM in 2011~\cite{GayzikEtAl11}, and have continued to develop models as shown in Tab.~\ref{tab_overviewTable}. 
The additions to the GHBMC family include both male and female models, child and adult models, occupant and pedestrian models, simplified and detailed models, and models of large and small people.
In addition, variants of these models have been developed such as the SAFER~HBM~\cite{PipkornEtAl21}, and the THUMS~TUC-VW AHBM~\cite{Yigit18}.
Even though there are many passive FE HBMs, these types of models continue to be developed and used during product development.

Current research focuses on making passive FE models that more accurately match population demographics, and can be used in a wider variety of simulation scenarios.
As the proportion of older adults in society increases, there is a focus on accurately capturing changes in age-related material properties~\cite{SchoellEtAl15} of the body.
There is also a growing focus on specialized algorithms to account for human variability by scaling existing models to represent the anthropometric properties of another person~\cite{LarssonEtAl19b, HwangEtAl16b}.
As the variety of simulation scenarios has grown, tools are being developed to position HBMs~\cite{KleinEtAl21, CostaEtAl21} in a specific posture, such as the PIPER framework~\cite{PIPER, KleinEtAl21} or the Oasys PRIMER positioning tool~\cite{CostaEtAl21}.
For some accident scenarios, however, the passive assumption of the human body is increasingly perceived as an unacceptable simplification.

\subsection{Active Finite Element Models}

Active safety systems begin to act in the phase prior to the crash, when accelerations are lower, and muscle activity influences the response of the body.
As summarized in~\cite{OesthLarssonJakobsson22} muscle activation influences the injury pattern in the pre-crash and in-crash phases:  neck muscle activity reduces the load on the cervical capsular ligament below the injury threshold and increases the load tolerance of the cervical spine from~1.8~kN to~3.1-3.7~kN; and active bracing of the leg's muscles increases the likelihood of injuries such as femoral shaft fracture.
Accordingly, efforts have been made to make active FE HBMs by including models of muscles and control systems.

Nearly all of the MBS and FE models mentioned in this work represent volumes of muscle by using arrays of Hill-type muscle cables that act independently.
While this approach reasonably approximates fiber paths in long skinny muscles, this approach cannot capture the curving paths that are produced when short bulky muscles bulge.
In contrast, the GHBMC model embeds Hill-type cables within a passive tissue volume making it possible to capture the fiber geometry even for many types of muscles.
Despite the existence of some advanced active 3D muscle models, none yet exist in production-ready code.

While muscle is mechanically complex, the way the nervous system controls muscle is perhaps even more challenging.
Since many active FE HBMs focus on a specific crash scenario, early models used open-loop control with pre-defined muscle activation patterns that were calculated offline from experimental electromyographic (EMG) data or by solving an application-specific optimization problem~\cite{IwamotoetAl12}. 
Closed-loop reflex control was later introduced to improve the response of the arms~\cite{OesthEtAl15a, Yigit18}, and the neck~\cite{OlafsdottirOesthBrolin17} during crash simulations.

Many additional control systems have been introduced to address the growing list of applications for active FE models.
The latest versions of the THUMS~5 and~6~\cite{IwamotoNakahira15, KatoEtAl18} contain simple~1D Hill-type muscles with angle-based controllers for posture maintenance. 
The VIVA+ model provides a posture controller for the neck that includes both vestibular and stretch reflexes~\cite{Kleinbach19, PutraEtAl19}. 
Similarly, the GHBMC models now include a number of stretch reflex controllers~\cite{FellerEtAl16, KleinbachEtAl17} and a head and neck controller that features both vestibular and stretch reflexes~\cite{CorreiaMclachlinCronin21}. 
The SAFER HBM~\cite{OesthLarssonJakobsson22}, based on THUMS~v3.0, has been tested with four different PID spinal posture controllers~\cite{OesthLarssonJakobsson22} in which a variety of pre-crash maneuvers and crash scenarios were evaluated. 
In addition, the SAFER~HBM now includes two error angular position feedback postural control schemes - one for the neck and one for the trunk.
Future research will likely continue to develop control systems to emulate in-vivo responses since muscle activity influences the pattern of injury~\cite{OesthLarssonJakobsson22}.

\begin{center}
  	\renewcommand*{\thefootnote}{\alph{footnote}}
	\begin{table}
	    \centering
		\begin{tabular}{@{}L{1.5cm}L{3cm}L{2.8cm}L{2.5cm}L{1.5cm}L{3.8cm}@{}} \toprule %
			\textbf{model} & \textbf{simulation times}\footnotemark[1] & \textbf{common analyses} & \textbf{domains} & \textbf{DOFs} & \textbf{variants}\footnotemark[2] \\ \midrule
			\textbf{RAMSIS} & $30\;\meh{s}$ to evaluate 50 positioning iterations & forward and inverse kinematics & ergonomics & $65$ & can be scaled to arbitrary anthropometries \\ %
			\multicolumn{6}{l}{
				\hspace{0.2cm}\makecell[l]{
					\begin{tabular}{p{2cm}p{13.2cm}}
						\emph{applications}: & evaluate anthropometrics-, posture-, visibility-, roominess-, discomfort-, reachability-, safety-, and operability-ergonomics criteria in the digital vehicle development~\cite{Wirsching19}
					\end{tabular}
				}
			}\\ \midrule

			\textbf{OpenSim} & - & inverse dynamics; optimal control & biomechanics; sports medicine & $21$ (Gait-Model) & can be scaled to arbitrary anthropometries \\ %
			\multicolumn{6}{l}{
				\hspace{0.2cm}\makecell[l]{
					\begin{tabular}{p{2cm}p{13.2cm}}
						\emph{applications}: & investigation of stiff-knee gait of a 12-year-old subject~\cite{DelpEtAl07} \\
						                     & inverse and forward dynamics and neuromuscular control of mammal motion~\cite{SethEtAl18} \\
						\emph{comments}:     & open-source
					\end{tabular}
				}
			}\\ \midrule

			\textbf{EMMA} & $45$ s to solve a simple optimal control problem & optimal control & early stage of safety system development & $40$ & Can be scaled to arbitrary anthropometries \\ %
			\multicolumn{6}{l}{
				\hspace{0.2cm}\makecell[l]{
					\begin{tabular}{p{2cm}p{13.2cm}}
						\emph{applications}: & simulation of dynamic driving maneuvers~\cite{RollerEtAl20}
					\end{tabular}
				}
			}\\ \midrule

			\textbf{AHBM} & $21$ min for $2$ s~\cite{KempterBechlerFehr20} & forward simulation & passive safety in impacts & $108$\footnotemark[3] & M50-O/P\\
			\multicolumn{6}{l}{
				\hspace{0.2cm}\makecell[l]{
					\begin{tabular}{p{2cm}p{13.2cm}}
						\emph{applications}: & occupant loading in rotated seats configurations~\cite{BeckerEtAl20}\\
						& Scenario bases safety performance assessment integrated vehicle safety systems~\cite{WimmerEtAl21}
					\end{tabular}
				}
			}\\ \midrule

			\textbf{GHBMC} & $96$ h for $304$ ms on $20$ CPUs\footnotemark[4] & forward simulation & injury prediction on tissue level & $2.3\; 10^6$ elem.& F05-O/P; M50-O/P; M95-O/P; \\ %
			\multicolumn{6}{l}{
				\hspace{0.2cm}\makecell[l]{
					\begin{tabular}{p{2cm}p{13.2cm}}
						\emph{applications}: & lumbar spine loads in reclined seating positions~\cite{RawskaEtAl19} \\
											& human response to underbody blasts~\cite{HostetlerEtAl21} \\
											& injuries of obese persons in sled tests an belt loading~\cite{GepnerEtAl18}\\
						\emph{comments}:     & commercial; free academic test-licenses
					\end{tabular}
				}
			}\\ \midrule

			\textbf{simpl. GHBMC} & $\sim 5$ times faster than full GHBMC~\cite{SchwartzEtAl15} & forward simulation & injury values on kinematic values & $0.8\; 10^6$ elem. & 6YO-P; F05-O/P; M50-O/P; M95-O/P;\\ %
			\multicolumn{6}{l}{
				\hspace{0.2cm}\makecell[l]{
					\begin{tabular}{p{2cm}p{13.2cm}}
						\emph{applications}: & feedback postural muscle control~\cite{CorreiaMclachlinCronin21} \\
						\emph{comments}:     & commercial; free academic test-licenses
					\end{tabular}
				}
			}\\ \midrule

			\textbf{THUMS} & $58$ h for $304$ ms on $20$ CPUs\footnotemark[4] & forward simulation & injury prediction on tissue level & $1.9\; 10^6$ elem. & 3YO-O/P; 6YO-O/P; 10YO-O/P; F05-O/P; M50-O/P; M95-O/P;\\ %
			\multicolumn{6}{l}{
				\hspace{0.2cm}\makecell[l]{
					\begin{tabular}{p{2cm}p{13.2cm}}
						\emph{applications}: & rib fracture mechanism of elderly occupant~\cite{TakahiraEtAl21} \\
											 & sport accidents involving children and adults~\cite{SmitEtAl21} \\
											 & influence of (autonomous) seating position, direction and angle on injuries~\cite{KitagawaEtAl17}\\
						\emph{comments}:     & open access since 2021; \emph{special~models}: elderly, obese, pregnant, whiplash
					\end{tabular}
				}
			}\\ \midrule

			\textbf{VIVA+} & $35$ h for $304$ ms on $20$ CPUs\footnotemark[4] & forward simulation & injury prediction on tissue level & $0.5\; 10^6$ elem. & 50F-O/P; 50M-O/P \\ %
			\multicolumn{6}{l}{
				\hspace{0.2cm}\makecell[l]{
					\begin{tabular}{p{2cm}p{13.2cm}}
						\emph{applications}: & addressing gender diversity, includes a detailed head-neck model ~\cite{OesthEtAl17a}\\
											& virtual testing in vehicle safety assessment of (standing) vulnerable road users~\cite{JohnEtAl22}\\
						\emph{comments}:     & open access, head and neck submodel for speedup available{\cite{PutraEtAl20}}
					\end{tabular}
				}
			}\\
			\bottomrule
			\multicolumn{6}{l}{
				\footnotemark[1]\footnotesize{unless otherwise specified, the values are given for calculations on a standard PC.}
			} \\
			\multicolumn{6}{l}{
				\footnotemark[2]\footnotesize{abbreviations: F: female; M: male; number corresponds to percentile; O: occupant model; P: pedestrian model;}
			} \\
			\multicolumn{6}{l}{
				\footnotemark[3]\footnotesize{$108$ rigid general coordinates, from 182 rigid bodies, 190 kinematic joints, 2174 FE elements from the skin}	
			}
		\\
		\multicolumn{6}{l}{
			\footnotemark[4]\footnotesize{further explanations on the calculation times of GHBMC, THUMS and VIVA+ are given in Sec.~\ref{subsec:case_study_safe_motorcycle}}	
		}
		\end{tabular}
		\caption{\enspace Overview of popular HBMs with which the authors have gained experience. The information in the table reflects the authors' view.\label{tab_overviewTable}}
	\end{table}
\end{center}

\section{Applications}\label{sec_application}

Digital HBMs are used in many technical and artistic fields.
To provide insight into the application of HBMs in the context of vehicle safety, the capabilities of some of the different HBMs, introduced in Sec.~\ref{sec:human_body_models_vehicle_safety}, are presented in three case studies.
In the first study, the application of a detailed FE model in the development phase of a new motorcycle safety system shows how complex models can contribute to a deep understanding of the events during a motorcycle accident and help in the evaluation of integrated safety systems.
The latter examples feature multibody systems used in the early development phase of new safety systems.
In a scenario-based assessment of the pre-crash phase, the Madymo AHBM is applied to gain a further understanding of the initial position of the head at the start of the in-crash phase.
Hereby, data from a Driver-in-the-Loop setup, showing the range and variance of motion in the pre-crash phase, is utilized to account for the occupants' behavior. 
In a third example, the EMMA model generates realistic movements with the help of optimal control to improve active and integral vehicle safety.

\subsection{Case Study -- Safe Motorcycle}\label{subsec:case_study_safe_motorcycle}

As omnidirectional tools, HBMs can be used in a wide range of applications. In the following, this is demonstrated with a comprehensive summary of the application of~state-of-the-art~FE~HBMs as motorcyclists in crash impacts from~\cite{MaierEtAl22}. Examples of current open and commercially available models are shown, a positioning workflow in the specific application for riders of a motorized two-wheeler is demonstrated and an impact response study on variations of HBMs, including a comparison versus computational models of ATDs, is discussed.

\paragraph{Positioning HBMs for Motorcycle Rider Postures}

\begin{figure}[tbp]
	\centerline{\includegraphics[width=17.5cm]{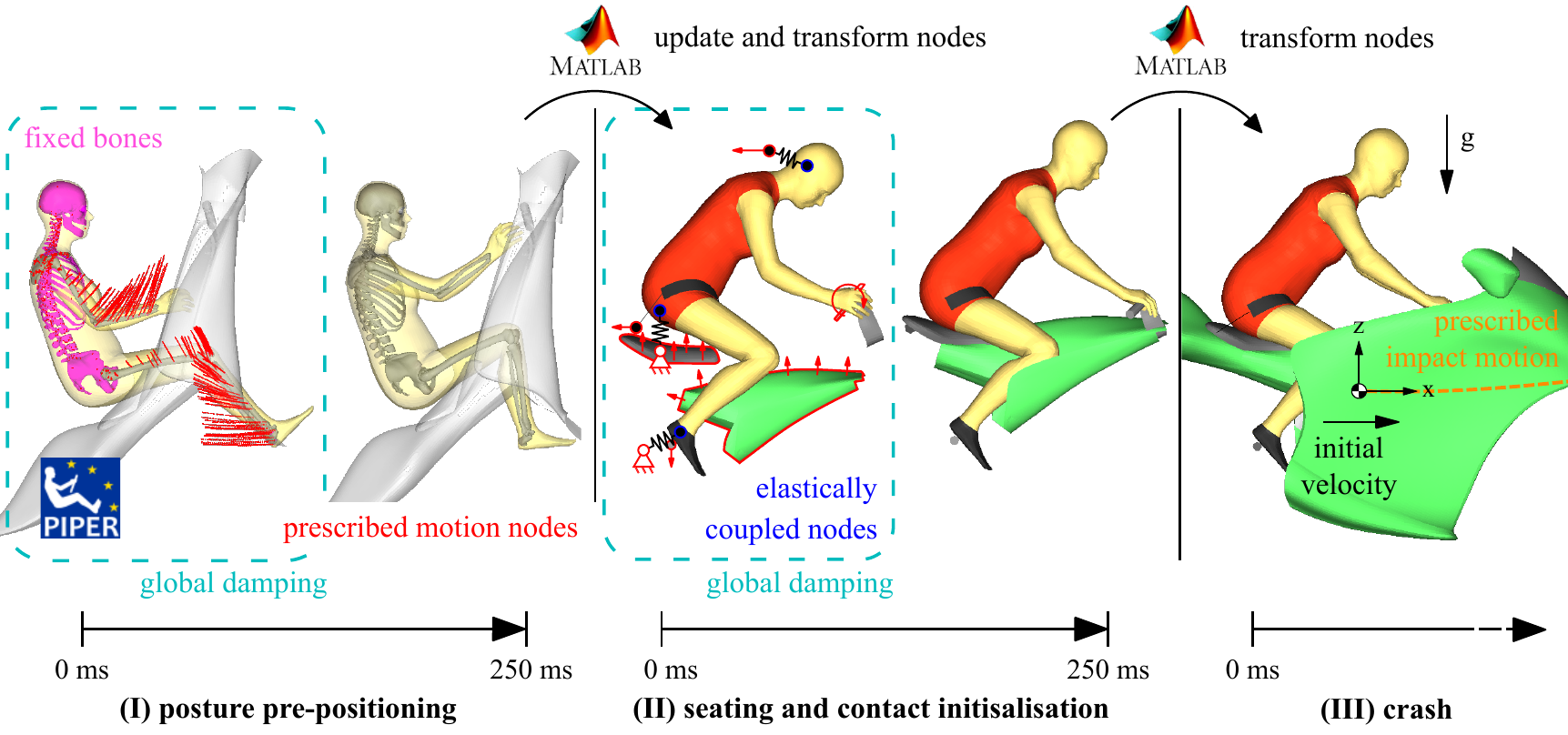}}
	\caption{\enspace Multi-step pre-positioning, seating, and contact initialization procedure for crash impact simulation of HBMs. \label{Fig:SafeMoto_HBMPositioning}}
\end{figure}

\begin{figure}
	\centerline{\includegraphics[width=17.5cm]{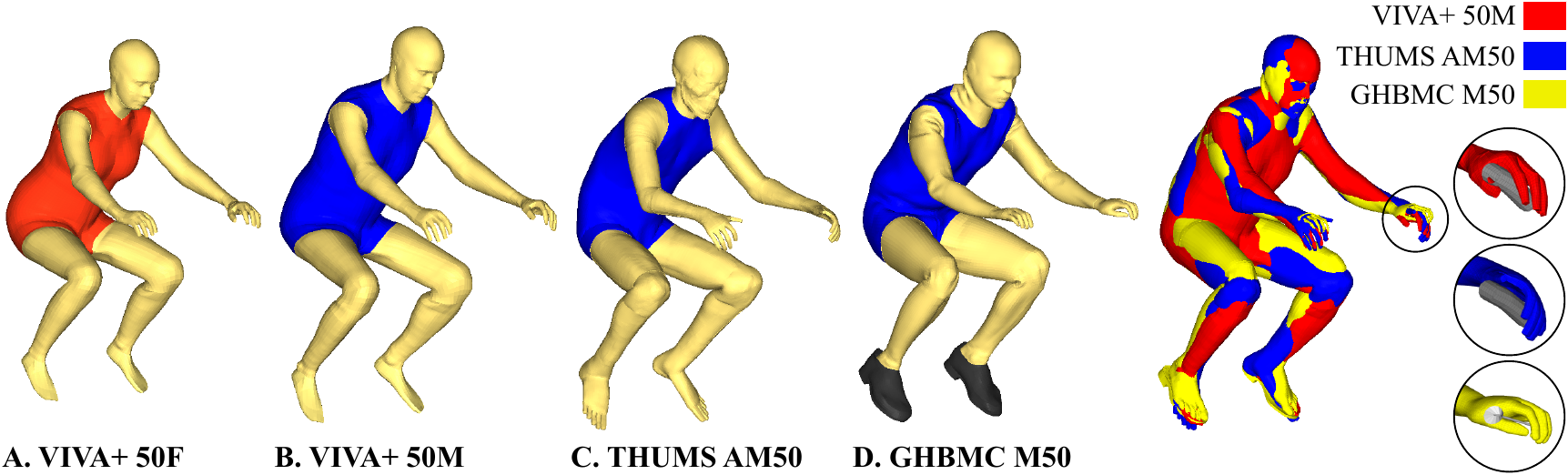}}
	\caption{\enspace Positioned and seated HBMs in the resulting motorcyclist’s riding postures~(\cite{MaierEtAl22}). \label{Fig:SafeMoto_HBMMotorcyclistPosition}}
\end{figure}

Positioning~HBM~FE models and initializing the actual simulation is an intricate and often laborious task. The default positions of most models are either standing postures of pedestrians or seated postures of car occupants. Three state-of-the-art models in the~LS-Dyna software environment are used: the~VIVA+~50F and~50M~\cite{JohnEtAl22,OpenVT21}, the THUMS~\cite{KatoEtAl18}, and the GHBMC~\cite{GHBMC22} models~(Fig.~\ref{Fig:SafeMoto_HBMDefault}). The posture of a motorcyclist is determined by their anthropometric measurements and the geometry of the motorcycle. The pelvis sits on the seat, the hands grip the handlebar, the feet rest on the footrests, and the head angle complies with the rider's line of vision. This posture is maintained by static muscle power. The developed workflow, which can be similarly applied to car occupants, consists of two preparation steps, pre-positioning~(I)~and seat and contact initialization~(II),~before the actual crash simulation~(III)~(Fig.~\ref{Fig:SafeMoto_HBMPositioning}).
\begin{itemize}
	\item \underline{(I)~Posture pre-positioning:} Each of the individual HBMs are prepositioned using the positioning and personalization tool PIPER~\cite{PIPER}. It is an academically developed open-source tool comparable to commercial solutions of other FE preprocessors based on HBM-specific metadata. For positioning, it provides functions for interactive manipulation of FE structures within the software through lightweight and, therefore, fast physics models and the definition of simulations for transient manipulation. Here, PIPER is used to define a transient simulation by prescribing the motion of the skeletal structure via elastic elements. The target position is determined by the joint angles of the lower and upper extremity joints with fixed landmarks for the spine and head. Assigning global damping in the explicit numerical simulation allows a faster positioning for a short simulation time.
	\item \underline{(II)~Seating and contact initialization:} For interaction with the motorcycle, a subsequent simulation is run. The rider-to-motorcycle contact surfaces, the seat, the handles, the footrests, and the other cockpit surfaces are initially scaled and distorted. Then, these geometries are morphed into the motorcycle’s actual geometry. This rotates the wrists and pushes the motorcycle against the human body. Prescribed spring elements are also used to push the head back, press the pelvis into the seat, and hold the feet to the footrests. Gravity is not applied yet.
	\item \underline{(III)~Crash:} The final step is the actual impact simulation. It starts shortly before impact deceleration, with an initial velocity and a prescribed impact trajectory for the motorcycle cockpit, similar to a crash pulse used for sled tests. The imposed motion is a multi-axial linear and angular motion. It is derived from a multibody simulation of the vehicle interaction of a motorcycle, a rider, and a passenger car as the accident opponent, explained in~\cite{MaierEtAl20}, as part of a simulation strategy that involves multiple models~\cite{MaierFehr23b}. The prescribed model is a motorcycle with additional passive safety equipment, airbags, thigh belts, and leg impact protection, as explained in~\cite{MaierEtAl21, MaierEtAl21b}. It restrains the rider to the two-wheeler, decelerates its movements, and prevents hard accident contacts by airbags.
\end{itemize}
Between the three separate simulation steps, the nodal coordinates of the mesh are transferred via Matlab procedures. This means mesh deformations are carried over into the subsequent simulation; other information, such as the element stresses due to deformations, is lost. The final postures at the end of step~(II)~agree reasonably well~(Fig.~\ref{Fig:SafeMoto_HBMMotorcyclistPosition}). For the female variant, adjusted footrests and handlebars are considered. The comparison demonstrates anthropometric differences of current~HBMs, although the resulting postures are similar overall. The comparison of the final hand positioning shows that the~HBMs differ in many modeling aspects. Here, the~VIVA+ and~GHBMC hands are non-deformable rigid bodies, which results in the fingers either piercing the handles, as shown for~VIVA+, or thinner handles must be formed, as applied for the~GHBMC. Only~THUMS allows for realistic gripping of the handles.

\paragraph{Variations of ATD and HBM Impact Response}

\begin{figure}
	\centerline{\includegraphics[width=17.6cm]{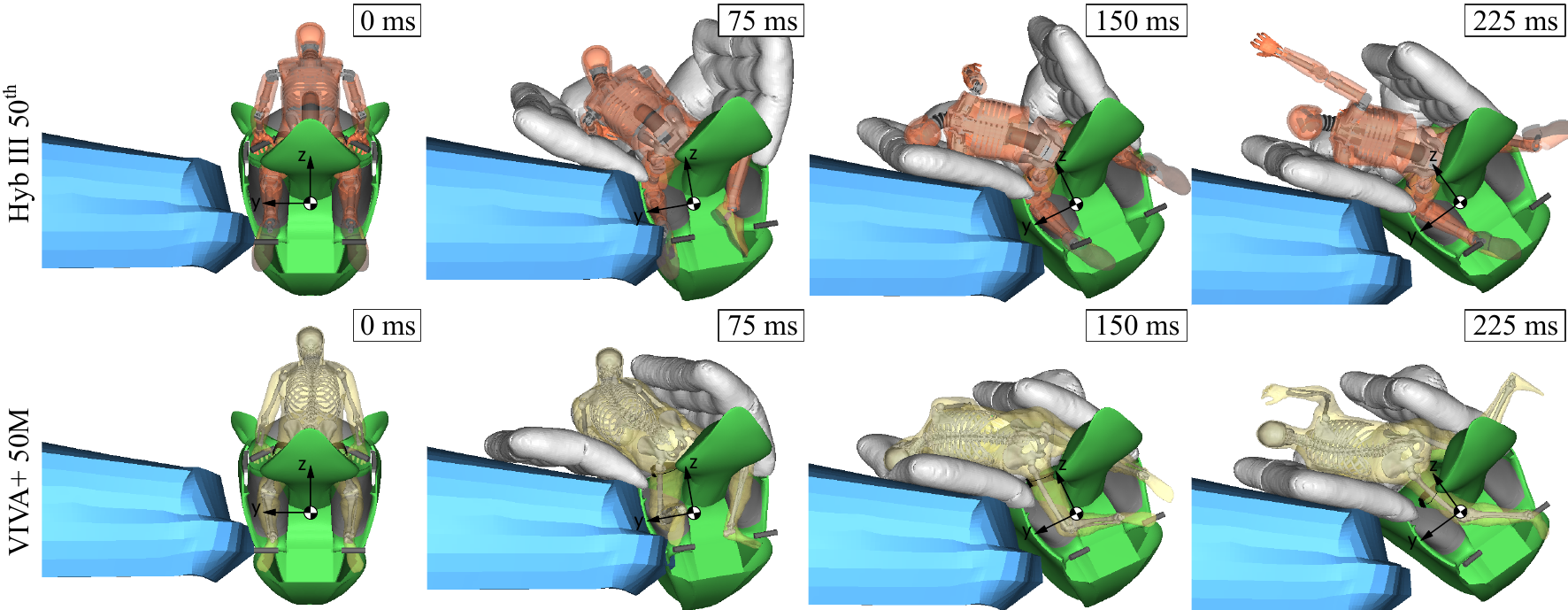}}
	\caption{\enspace Rider trajectories of Hybrid~III 50\textsuperscript{th} percentile~ATD~(top) and VIVA+~50M~HBM~(bottom) for lateral car impact~(\cite{MaierEtAl22}). \label{Fig:SafeMoto_ATDvsHBMImpactComparison}}
\end{figure}

\begin{figure}
	\centerline{\includegraphics[width=17.6cm]{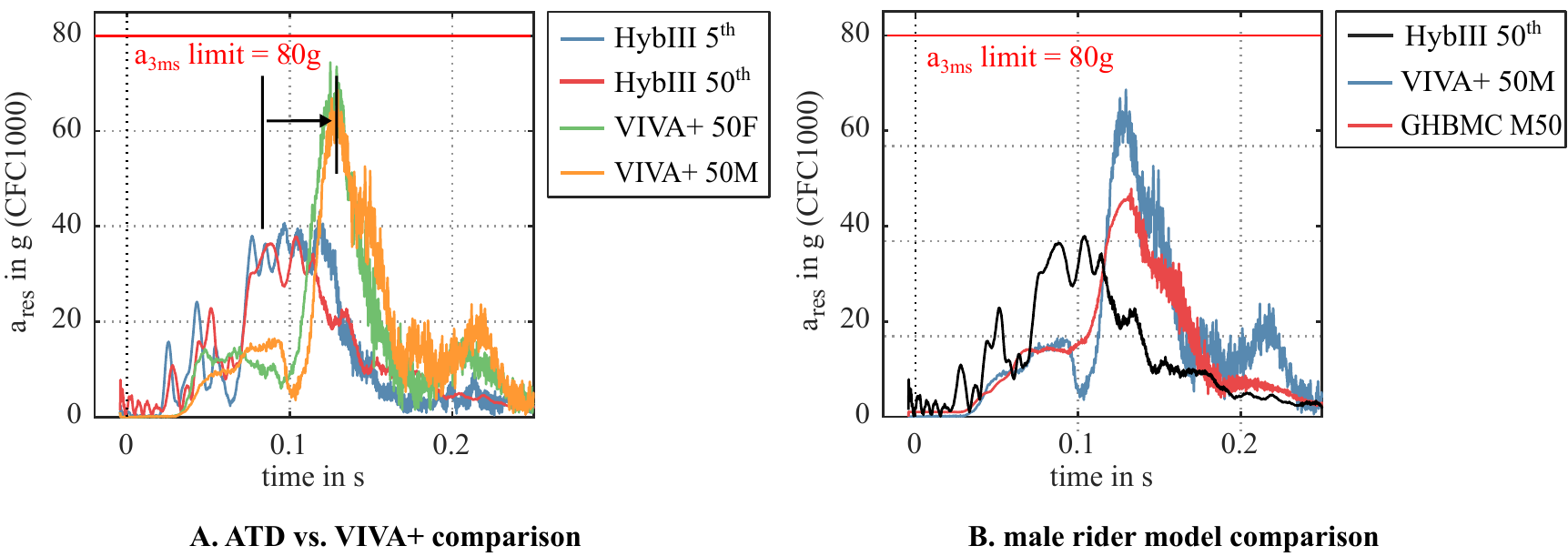}}
	\caption{\enspace Resulting head~COG linear accelerations for ATDs and HBMs. Note: THUMS~AM50 does not provide a head~COG acceleration sensor~(\cite{MaierEtAl22}). \label{Fig:SafeMoto_LinAccComparison}}
\end{figure}

Evaluating the primary impact kinematic response demonstrates significant differences between HBMs and ATDs in a lateral impact of a passenger car with~35~km/h against the stationary motorcycle~(Fig.~\ref{Fig:SafeMoto_ATDvsHBMImpactComparison}). The models are partially displayed transparently to illustrate the mechanical and human skeletal structure. The side airbag deploys within~35~ms and is pushed down by the riders, with the car acting as a reaction surface. The spine of the HBM is laterally more flexible than the mechanical replication of the ATD, observable at $75\;\meh{ms}$ by a delayed head motion and at~150~ms by a more continuously deflected cervical and thoracic spine. The head is accelerated from the upper body like a whip. Virtual sensor histories~(Fig.~\ref{Fig:SafeMoto_LinAccComparison}~(A)), reveal higher linear accelerations at a later time. The result shown is expected, as~ATD are designed to be biofidelic in very specific load cases, mostly to depict car occupants. Here, human models are a meaningful omnidirectional supplement.

\begin{figure}
	\centerline{\includegraphics[width=17.6cm]{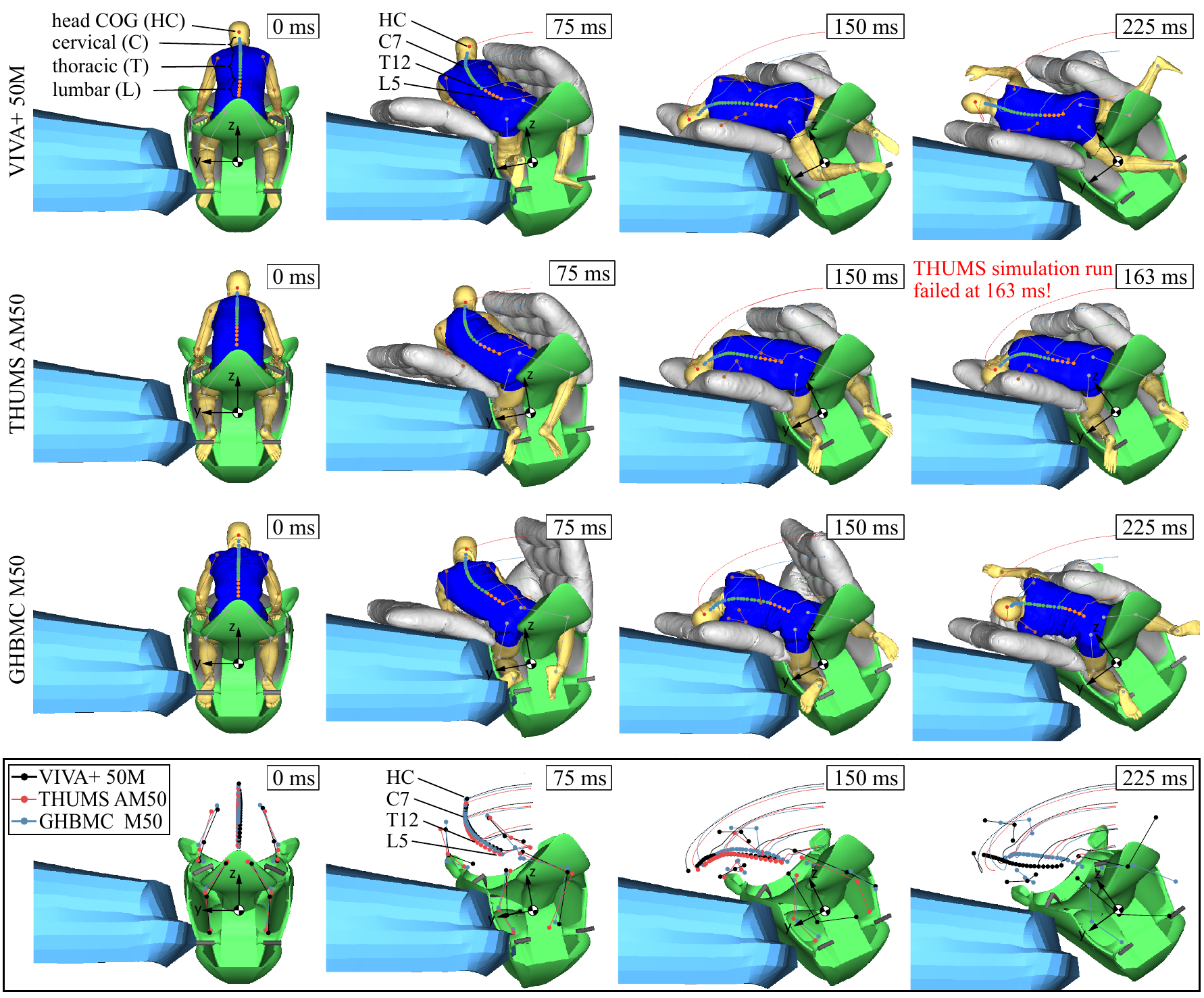}}
	\caption{\enspace Rider trajectories~(top) and comparison~(bottom) of male HBMs for lateral car impact~(\cite{MaierEtAl22}). \label{Fig:SafeMoto_MaleHBMImpactComparison}}
\end{figure}

Comparing multiple HBM variants shows minor variations in impact kinematics~(Fig.~\ref{Fig:SafeMoto_MaleHBMImpactComparison}). At each snapshot, the main skeletal structure is visible through the semi-transparent body. The trajectories of head~COG, C7, T12, and L5 are continuously traced. Overall the skeletal structures show similar initial position and impact responses for the~HBMs upon impacting the airbags. Differences arise from variations in anthropometry and mass distribution from underlying body scans and/or the modeling of materials, joints, and contact definition. For linear acceleration, this leads to some greater variations~(Fig.~\ref{Fig:SafeMoto_LinAccComparison}~(B)). Like different real-world humans with similar - but of course not identical - physiques, this results in a corridor that, in turn, challenges the robustness of a safety solution. Working with several models to represent expected variance may be advised.

The available models vary significantly in the calculation times, see Tab.~\ref{tab:Calctime_safe_Motorcycle}. The VIVA+ model with only 490K elements is the computational least costly model of the three. The VIVA+ model is not as densely meshed as the THUMS and the GHBMC model. The GHBMC model with 1.9M nodes is the most complex. The calculation factor increases to over 1\,000\,000 for the crash simulation, which means that for $1$ second simulation time over $1\,000\,000$ seconds elapse until we receive the results. For all models, the calculation factor is smaller in the pre-positioning and the seating and contact initialization problems.

\renewcommand\arraystretch{1.2}
\begin{table}[!ht]
	\caption{\enspace Calculation times for the safe motorcycle study with HBMs. All simulations where performed with LS-Dyna R 9.3.1 mpp single precision version using 20 CPUs.\label{tab:Calctime_safe_Motorcycle}}
	\centering
	\begin{tabular}{@{}lrr@{}}
	\multicolumn{2}{l}{(I) prepositioning (problem time = $250\;\meh{ms}$)} & calc. factor \\ \hline
		VIVA+ 50M				& $7.8\;\meh{h}$  		& $\sim$\,113\,000 \\ 
		THUMS AM50			&$19.6\;\meh{h}$  		& $\sim$\,282\,000 \\ 
		GHBMC M50-O 		    &$57.0\;\meh{h}$  		& $\sim$\,821\,000 \\ 
		\multicolumn{3}{l}{(II) seating and contact initialization (problem time = $250\;\meh{ms}$)} \\ \hline
		VIVA+ 50M				& $13.2\;\meh{h}$ 		& $\sim$\,190\,000 \\ 
		THUMS AM50			& $20.6\;\meh{h}$ 		& $\sim$\,296\,000 \\ 
		GHBMC M50-O 		    & $60.5\;\meh{h}$	    & $\sim$\,872\,000 \\ 
		\multicolumn{3}{l}{(III) crash simulation (problem time = $304\;\meh{ms}$, scenario 7)} \\ \hline
		VIVA+ 50M				& $35.2\;\meh{h}$ 		& $\sim$\,420\,000 \\ 
		THUMS AM50			& $58.2\;\meh{h}$		& $\sim$\,690\,000 \\ 
		GHBMC M50-O 		    & $95.5\;\meh{h}$		& $\sim$\,1\,130\,000 \\ 
	\end{tabular}
\end{table}

\begin{figure}
	\centerline{\includegraphics[width=17.6cm]{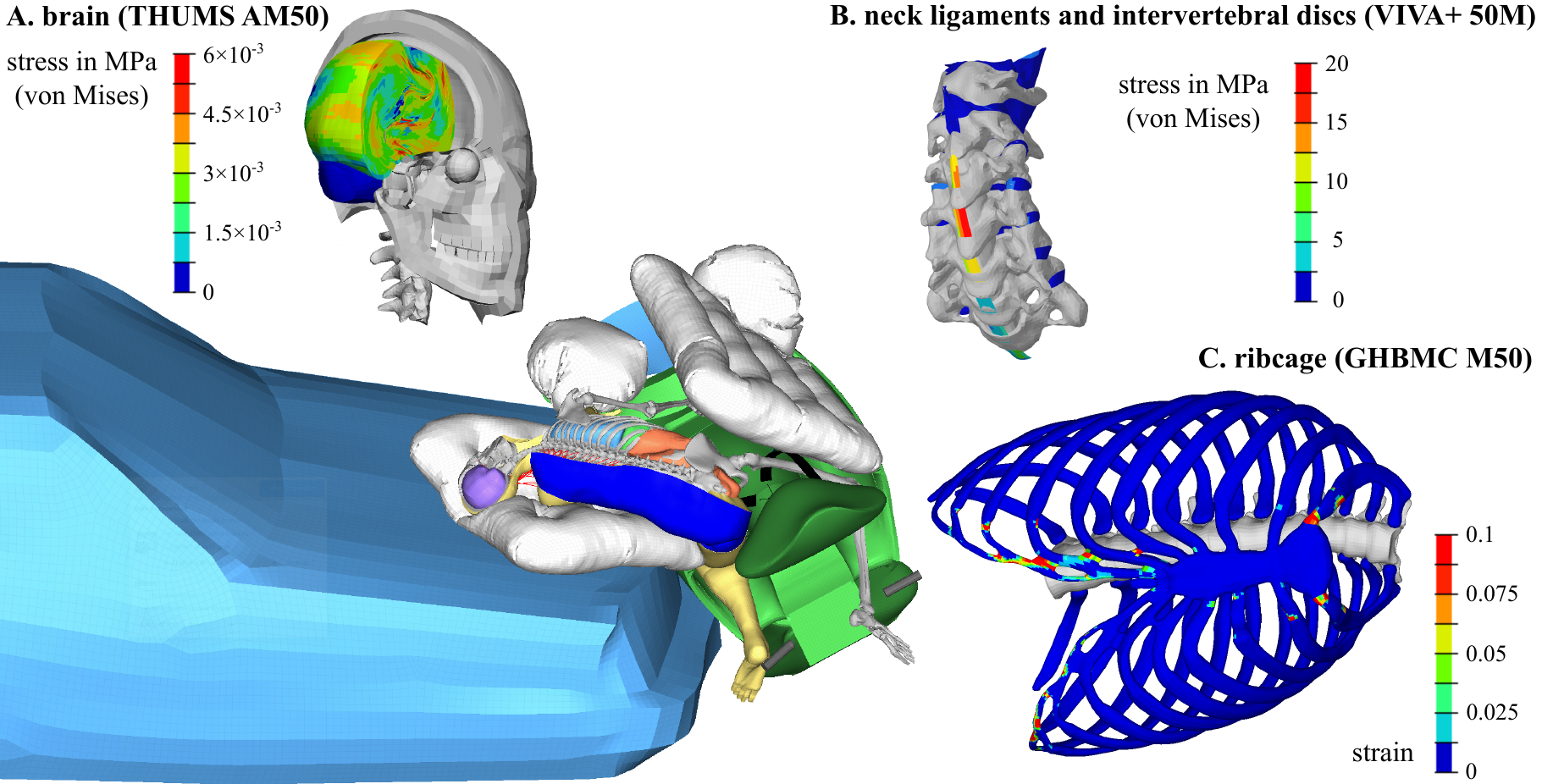}}
	\caption{\enspace THUMS accident response~(partially blanked) with illustrations of exemplary strain- and stress-based assessments of various body parts of different~HBMs. \label{Fig:SafeMoto_VirtualSensor}}
\end{figure}

HBMs enable a close examination of injury mechanisms. Typically, injury criteria like the head injury criterion~(HIC)~\cite{KleinbergerEtAl98} are derived from sensor loads. These criteria are associated with injury probabilities such as from~\cite{NHTSA95}. These, in turn, are based on injury classifications like AIS~\cite{AAAM08}. Based on these injury probabilities, limits are set for recommendations of consumer ratings or regulations to certify passive safety systems. Here, HBMs provide much more profound insight. They offer the possibility of a virtual magnifying glass that allows observing human injury mechanisms that are strain or stress-based~(Fig.~\ref{Fig:SafeMoto_VirtualSensor}).

\subsection{Case Study -- Madymo AHBM in a Driver-in-the-Loop Simulator and Scenario-based Assessment\label{subsec:case_study_safety_simulation_framework}}

People adopt a variety of postures~\cite{FiceBlouinSiegmund18} during naturalistic driving. 
The differences in initial posture in the pre-crash phase can have an influence on the risk of injury during the in-crash phase~\cite{OesthLarssonJakobsson22}.
Muscle activity strongly influences posture during the pre-crash phase because the accelerations are considerably lower when compared to the in-crash phase (Fig.~\ref{fig_DrivingSituations}).
Control systems that emulate the reflexes measured during in-vivo driving experiments have been developed so that the resulting changes in posture can be predicted~\cite{PutraEtAl20,KempterEtAl20}.

Tuning postural and reflex controllers~\cite{PutraEtAl20,KempterEtAl20} can be challenging because of the large number of simulations required to converge on a solution.
Although surrogate models~\cite{HayEtAl22, Hay22} are being developed to reduce simulation times further, it is still typical to simulate both multibody and finite element models using numerical integration.
The influence of reflexes on the initial posture of the occupant can be easily studied by simulating a multibody model because the calculation factor is so much lower than a finite element model.
A popular choice to simulate the pre-crash phase is Simcenter Madymo AHBM which can simulate $2\;\meh{s}$ of problem time in $21\;\meh{min}\;4\;\meh{s}$ of calculation time ~\cite{KempterBechlerFehr20} \footnote{Madymo R7.7 on a 12 core Intel i7-8700 processor with $3.20\;\meh{GHz}$ and $64\;\meh{GB}$ RAM} resulting in a calculation factor of 632 which is far less than the calculation factor of 100,000-1,000,000 that appeared in the finite element models used in Sec. \ref{subsec:case_study_safe_motorcycle}.
While tuning a postural or reflex controller is challenging, it is also difficult to obtain naturalistic kinematic and EMG responses during the pre-crash phase.

\begin{figure}[tbp]
	\begin{minipage}[b]{0.52\textwidth}
		\centering
		\includegraphics[width=1.0\textwidth]{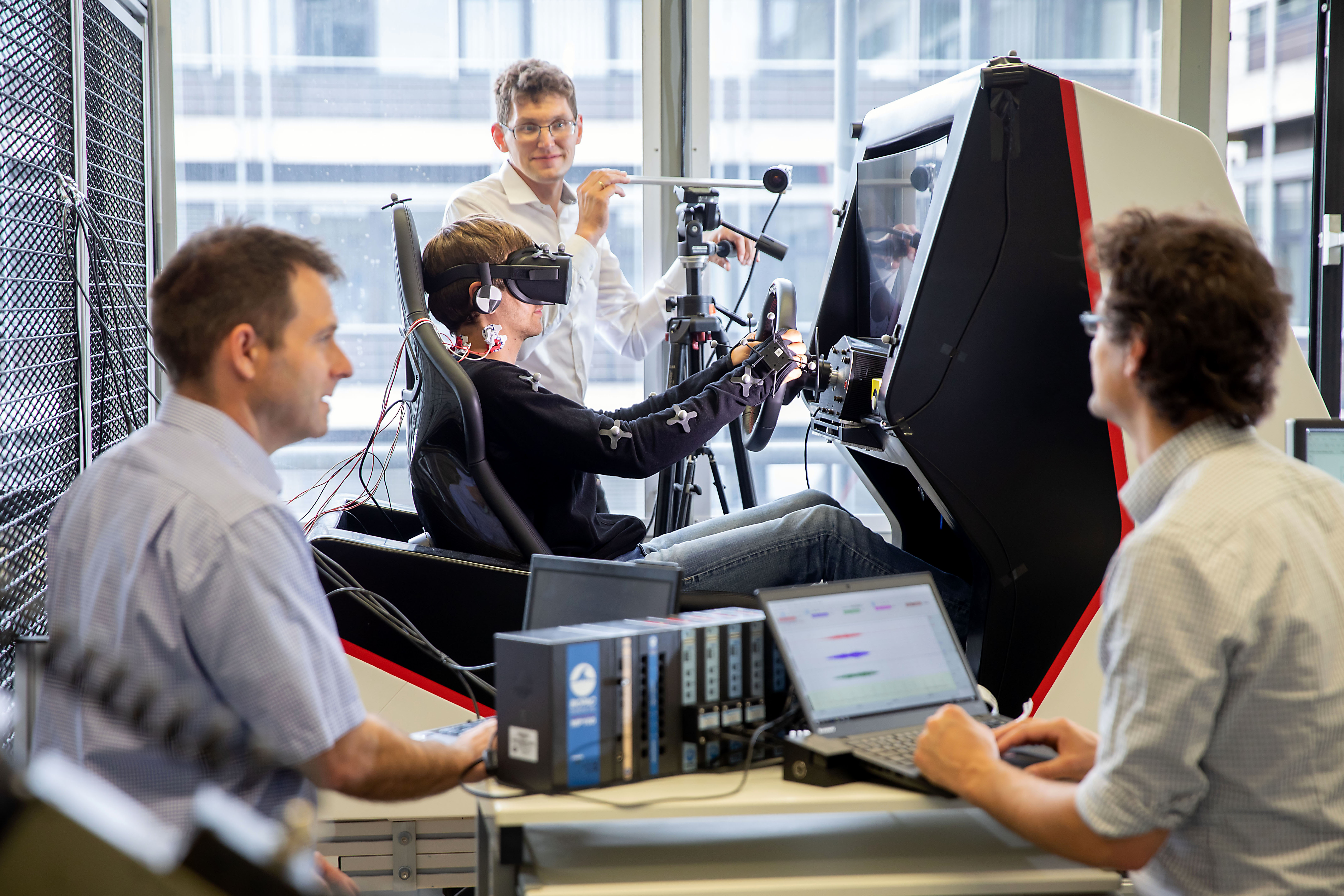}
	\end{minipage}
	\hfill
	\begin{minipage}[b]{0.46\textwidth}
		\centering
	    \includegraphics[width=1.0\textwidth]{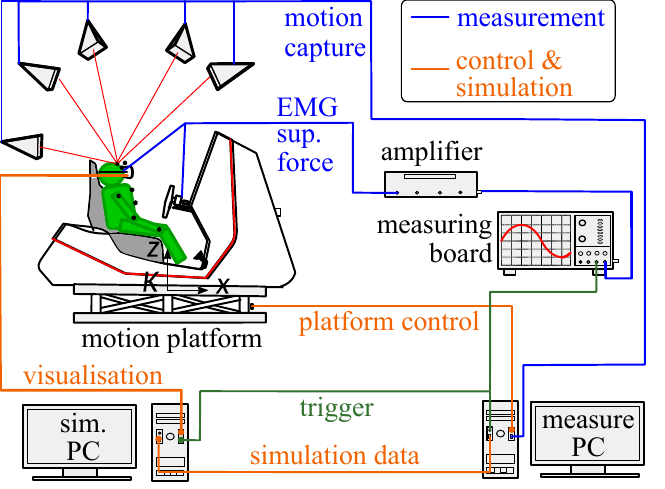}
	\end{minipage}\\
	\begin{minipage}[t]{0.52\textwidth}
		\centering
		\textbf{A. Driver-in-the-Loop (DiL) simulator setup}
 	\end{minipage}
	\hfill
	\begin{minipage}[t]{0.46\textwidth}
		\centering
		\textbf{B. Internal and external signal flow for simulation and experimental data acquisition}
	\end{minipage}\\
	\caption{\enspace Setup of the Driver-in-the-Loop (A) simulator with a Porsche racing simulator, a Stewart motion platform and an Optitrack motion capture system consisting of four cameras attached to the ceiling and EMG measurement where the electrodes are placed on the neck muscle \copyright University of Stuttgart/Uli Regenscheid. Internal and external signal flow for simulation and experimental data acquisition (B) to collect synchronized measurements (blue) of the human driver (green) sitting inside the simulator.\label{fig_DIL_Fehr}}
\end{figure}

\begin{figure}[tbh]
	\centerline{\includegraphics[scale=0.65]{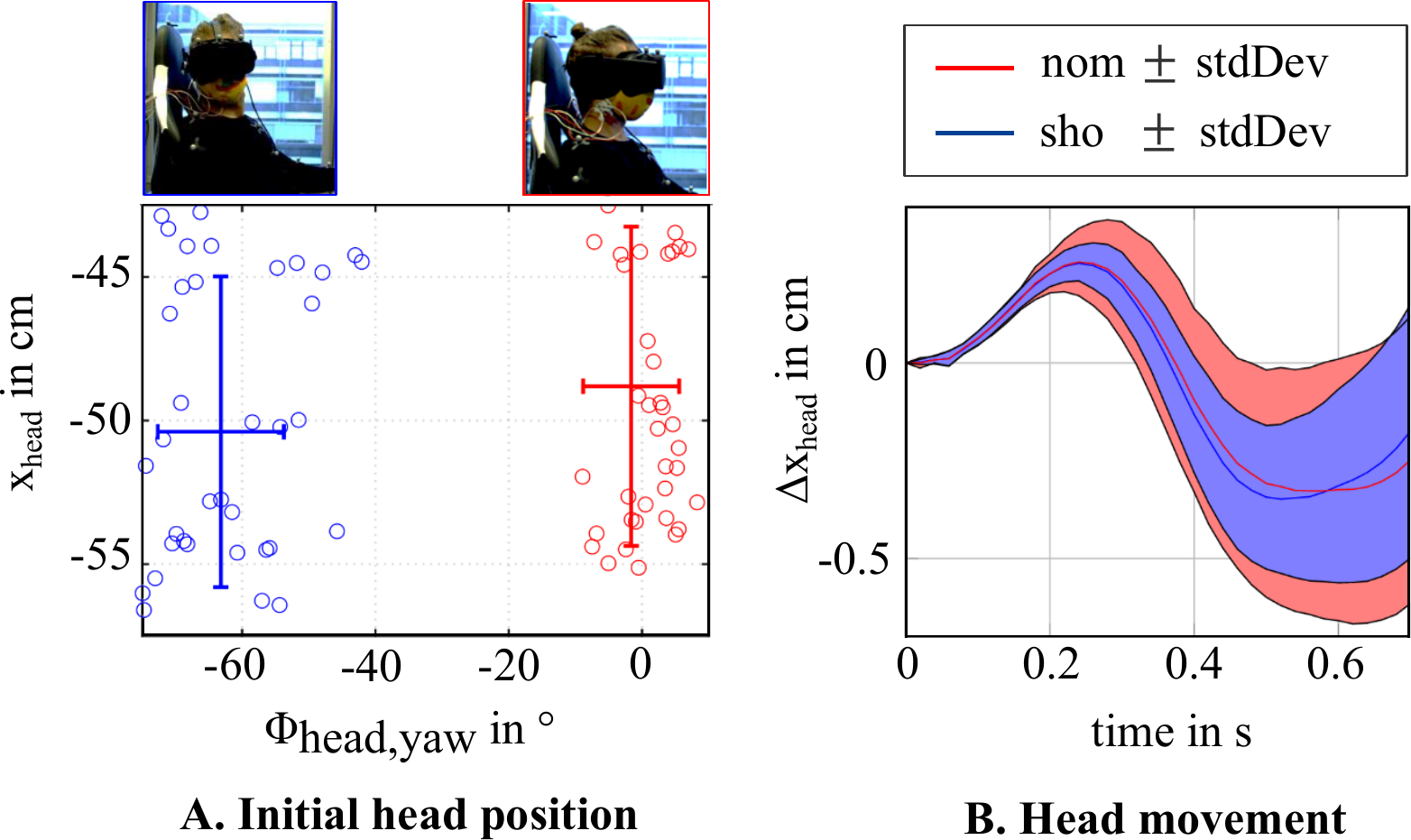}}
	\caption{\enspace Initial head postures (A) shown as initial yaw angle and longitudinal position for nominal view (red) and shoulder view (blue) as mean $\pm$ standard deviation.
	Head longitudinal displacement over time during braking (B).\label{fig:Head_Postures_Variances}}
\end{figure}

A Driver-in-the-Loop (DiL) simulator ~\cite{KempterBechlerFehr20, KempterEtAl23, FehrEtAl21b} makes it possible to systematically study the effects of posture~\cite{KempterEtAl21} during the pre-crash phase, while making laboratory-quality measurements of participant kinematics and muscle activity.
In a volunteer test study approved by the local ethical committee, the behavior of 17 typical participants was investigated during a mixed reality braking scenario. 
The participants were asked to look ahead (nom), or to the right (sho) (see photos in Fig. \ref{fig:Head_Postures_Variances} A.), so that their movements and muscle activity could be studied~\cite{KempterEtAl18, KempterEtAl23}. 
As in other volunteer studies, a large variability in the volunteers’ behavior can be observed (Fig. \ref{fig:Head_Postures_Variances} A. and B.).

Madymo models have been developed that can reproduce a similar degree of head-neck movement variability as is observed during in-vivo experiments.
In these simulations, the platform motion is constrained to follow the measured movements of the DiL cockpit.
A controller is formulated that uses plausible sensory information, such as muscle stretch \cite{KempterEtAl23}, which drives a parametric controller to activate the muscles of the neck.
With a parametric space defined, optimization is used to search for the parameter values that minimize the differences in simulated and experimental kinematics, as well as EMG onset times.
The individual subject-specific parameters can be used to define parameter ranges to replicate the variability of the human behavior observed in the experiments. 
Instead of one nominal head position during a pre-crash maneuver, the simulation predicts a range of possible head postures accounting for the identified variability in the occupants' behavior. 
This could enable the development of more advanced safety systems that function effectively across a range of occupant postures.

Scenario based testing considers the variability in acceleration introduced by the numerous trajectories that a vehicle could adopt within a specific driving scenario.
Ideally both active and passive safety systems are included during scenario-based testing because these systems can influence the trajectory and acceleration of the vehicle.
A simulation framework~\cite{HayEtAl22, WimmerEtAl21} has been developed to systematically define a range of scenarios~\cite{HayEtAl22, WimmerEtAl21} that includes an occupant pre-crash and in-crash simulation. By including both a pre-crash and in-crash simulation it is possible to trace the conditions of the pre-crash phase to the injury outcomes during the in-crash phase.

\begin{figure}[tbp]
	\centerline{\includegraphics[scale=1.0]{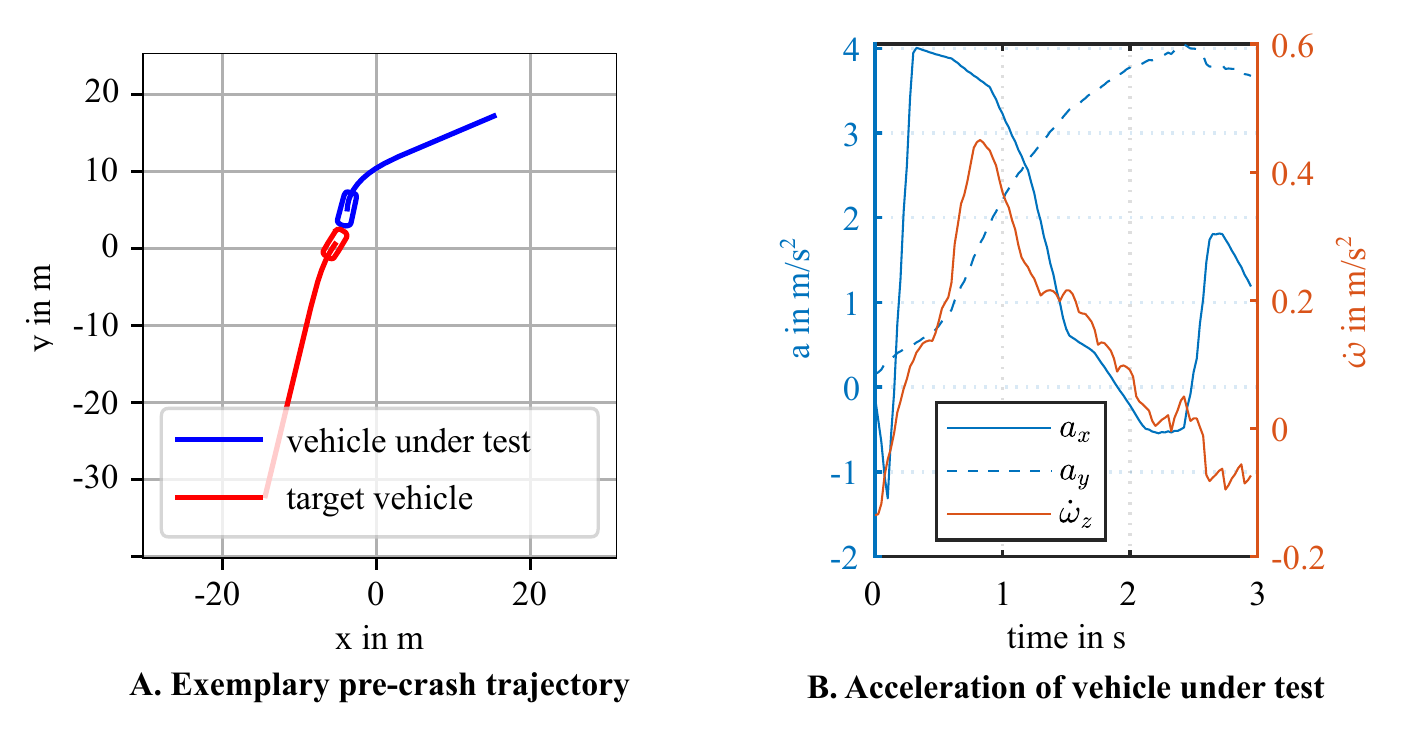}}
	\caption{\enspace (A) Exemplary pre-crash trajectory of the vehicle under test (blue) and the target vehicle (red); (B) Corresponding $x$-, $y$-, and $\omega_z$- acceleration of the vehicle under test.\label{fig_Madymo_one_scenario}}
\end{figure}

Vehicle scenarios can be defined and simulated using the commercial software tool ANDATA (Fig.~\ref{fig_Madymo_one_scenario} A.).
These trajectories can be refined by using an advanced vehicle model, such as one from Carmaker, to try to follow the trajectory from ANDATA. 
This refinement is important so that the effects of autonomous emergency braking (AEB), and other safety systems, can be included to make the trajectory more realistic.
In addition to improving the realism of the trajectory, the Carmaker simulation can be used to determine if a crash is unavoidable, or can be prevented, by an active safety system.
The refined trajectory (Fig.~\ref{fig_Madymo_one_scenario} B.) can be used as input for an occupant simulation.

\begin{figure}[tbp]
	\centerline{\includegraphics[scale=1.0]{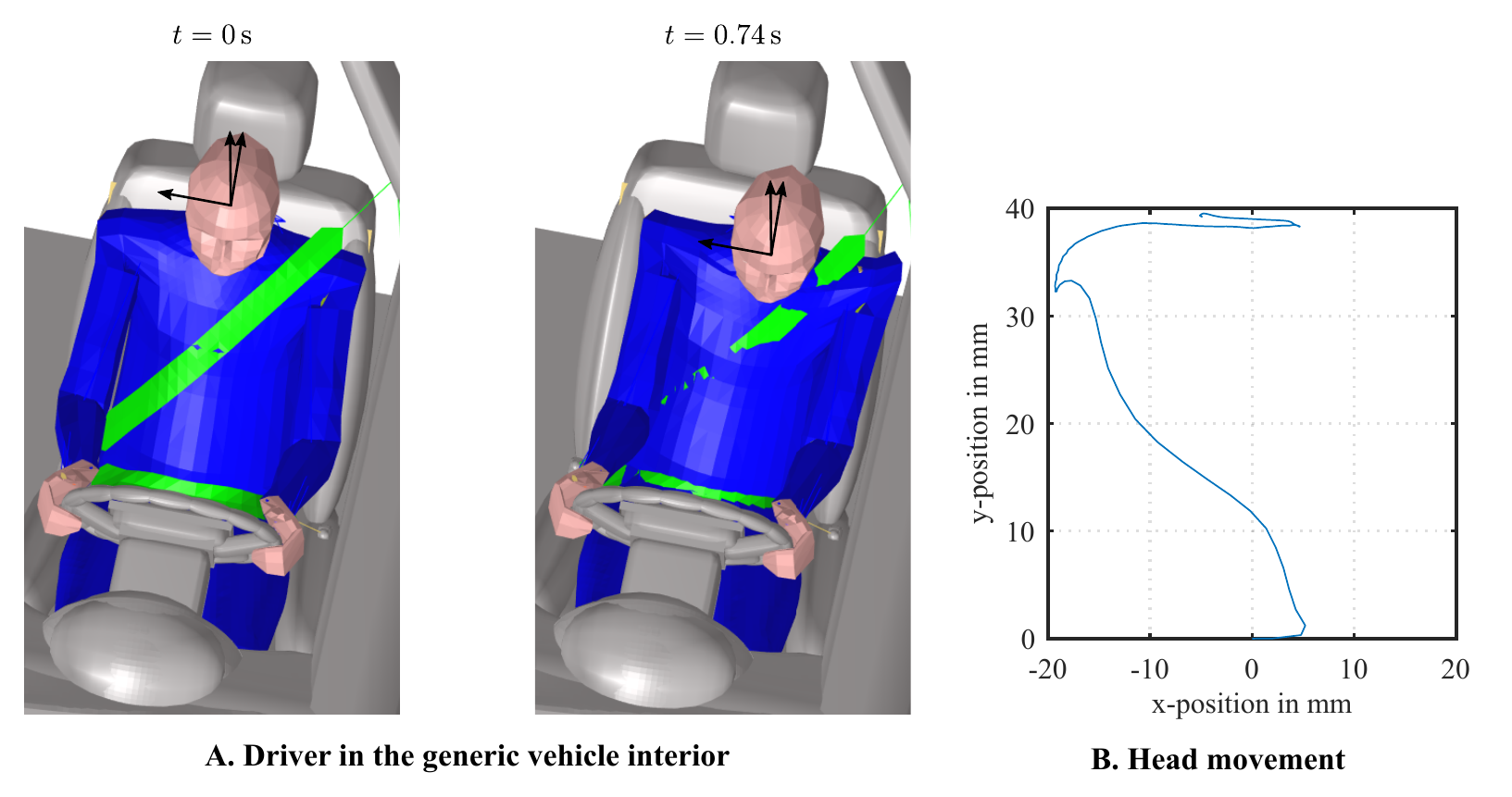}}
	\caption{\enspace (A)~Two time frames of the simulation of the driver in the generic vehicle interior with fastened seat belt; (B)~Head movement during simulation\label{fig_Sim_Madymo_PreCrash_Simualtion}}	
\end{figure}

The occupant simulation is used to determine how pre-crash movements affect the risk of injury during the later in-crash simulation.
The occupant simulation makes it possible to include a realistic model of the vehicle interior~\cite{BeckerEtAl20} (Fig.~\ref{fig_Sim_Madymo_PreCrash_Simualtion}) including the control panel, steering panel, and driver seat.
By using multibody models for the simulation, it makes possible to simulate long pre-crash phases (Fig.~\ref{fig_Sim_Madymo_PreCrash_Simualtion} shows the occupant response to the $740\;\meh{ms}$ long pre-crash phase illustrated in Fig.~\ref{fig_Madymo_one_scenario}) that would be prohibitively expensive to evaluate using an FE model. 
The final position of the pre-crash phase is then the input to the in-crash simulation, which is then made either with FE models or the AHBM.
The AHBM can also be used in the in-crash phase to predict injury values because it has been validated against PMHS experimental data.

\begin{figure}[tbp]
	\centerline{\includegraphics[scale=1.0]{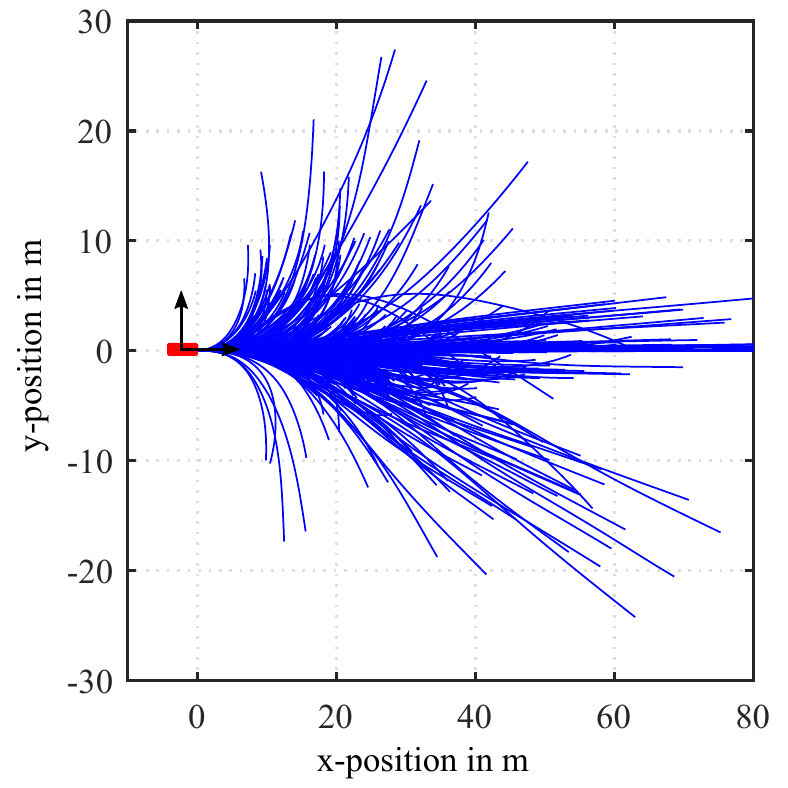}}
	\caption{\enspace Considered pre-crash trajectories of the vehicle under test in blue, resulting from the scenario simulation in the event of a frontal collision with the target vehicle in red. The vehicle under test trajectories (in blue) are illustrated relative to the position of the target vehicle (red). Only those pre-crash trajectories resulting in a frontal collision are shown.
		\label{fig_all_impact_scenarios}}
\end{figure}

\begin{figure}[h]
	\centerline{\includegraphics[scale=1.0]{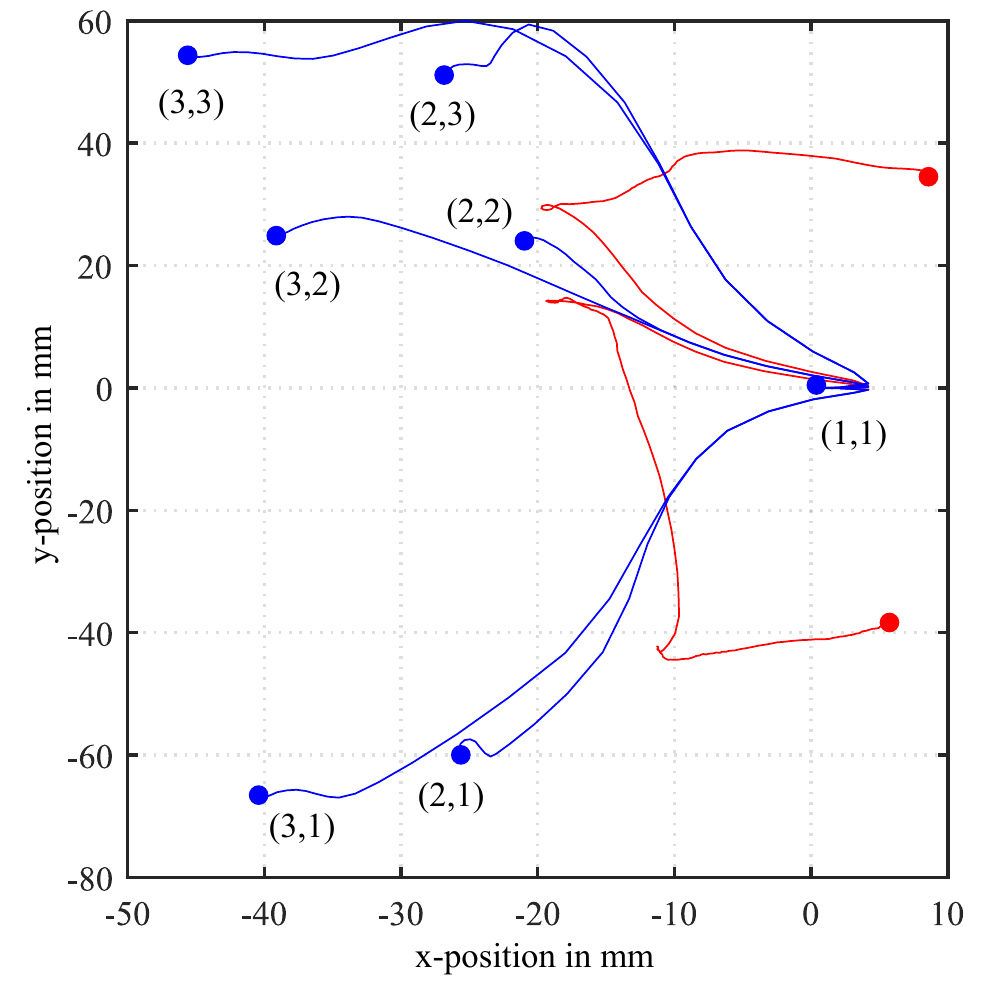}}
	\caption{\enspace Head displacements describing the occupants' range of movement resulting from simulations with generic vehicle acceleration signals in blue with corresponding indices i, j and driving scenarios with positive longitudinal acceleration in red. The end position at the time of collision is indicated by the dots.\label{fig_Range_of_Head_Motion}}	
\end{figure}

It remains challenging to simulate a typical range of scenarios (Fig.~\ref{fig_all_impact_scenarios}) with the AHBM and evaluate the occupant's positions.
As an example, there are a number of combinations of generic pulses that might be simulated to represent accelerations during AEB braking through a curve:
(i)~the $C0$ steady lateral x-acceleration in the tuple with $a_x^{max}=\{0,\, 4.25,\,8.5\} \frac{m}{s^2}$ ramped up with $\dot{a} =23\, \frac{m}{s^3})$, 
(ii)~the $C0$ steady transversal y-acceleration in the tuple with $a_y^{max}=\{-4.5,\, 0,\, 4.5\}\frac{m}{s^2}$ with a maximum acceleration in the y-direction of $4.5\, \frac{m}{s^2}$ ramped up with $\dot{a} =23\, \frac{m}{s^3})$ (iii)~combination of both; and (iv) scenarios form the scenario pool (see Fig.~\ref{fig_all_impact_scenarios}) to represent a maximum positive lateral acceleration plus mirrored transversal accelerations.
For each combination of the generic acceleration signals in x- and y-direction with maximum $a_{i,j}^{max}=( a_{x,i}^{max}, a_{j}^{max})$ a Madymo simulation is performed with $t \in [0 , 3 ]\;\meh{s}$.
It can be seen that negative transversal acceleration leads to a negative y-position (head moves to the left, Fig.~\ref{fig_Range_of_Head_Motion}).
The AHBM tends to move his head slightly to the right (positive head position).
Furthermore, the head moves to the front with breaking.
A clear correlation is observable which is a clear indication that passive and active safety devices should be considered holistically, as one influences another.
Overall the controller of the AHBM can minimize the head and occupant motion to a realistic small movement.
Without the stabilization provided by the muscles, the head would move a lot.

\subsection{Case Study -- \texorpdfstring{EMMA}{EMMA}: Investigation of the Reachability of the Steering Wheel from a Reclined Posture}\label{sec:case_study_EMMA4Drive}

The following exemplary utilization of the EMMA model investigates the reachability of the steering wheel from a reclined posture with and without actively moving the backrest to support the passenger.
The reclined starting pose imitates a passenger in a zero-gravity posture with the arms laying on armrests.
This position is assumed to become common in highly automated vehicles when the occupant is not actively controlling the vehicle~\cite{OestlingLarsson19}.
Nevertheless, for level three and four automated driving, the occupant must be able to take over control of the vehicle in case of an emergency~\cite{NHTSA13}.
To evaluate the integration of an actively moving backrest, the time to straighten up and reach the steering wheel is examined with and without support from the backrest.

\begin{figure}[t]
	\centerline{\includegraphics[width=0.3\linewidth]{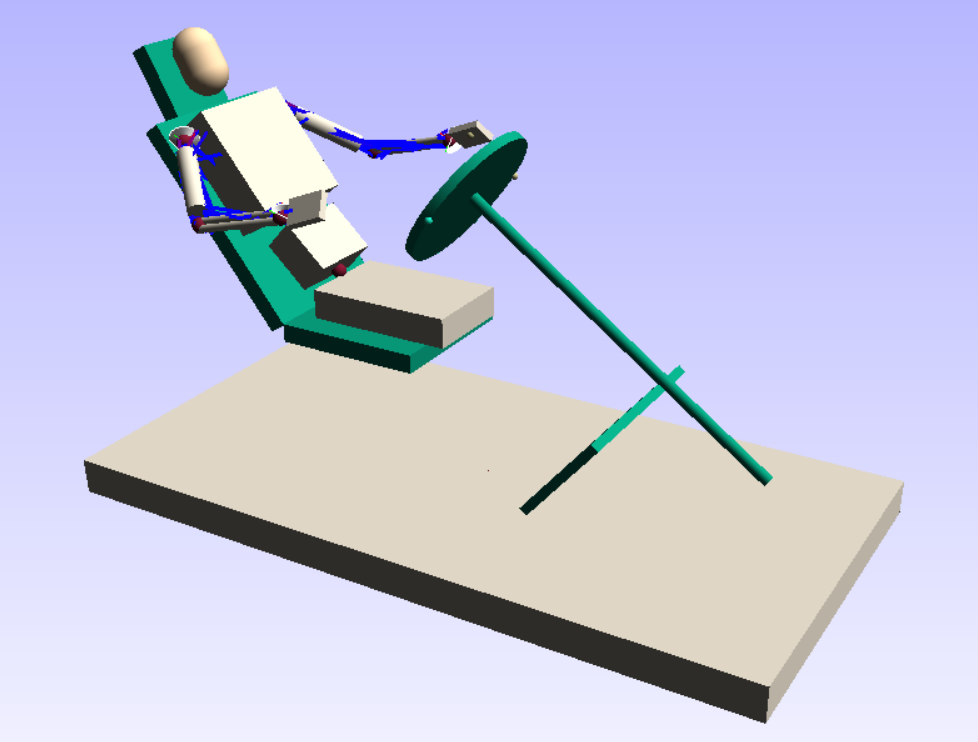}}
	\caption{\enspace Visualization of the EMMA model in the initial reclined posture that is used for the investigation of the reachability of the steering wheel.\label{fig_E4D_initial}}
\end{figure}

The EMMA model is used to simulate the movement of the upper body of the occupant from the reclined posture to the steering wheel.
Therefore, a model of the upper body is used and -- as shown in the initial reclined posture in Fig.~\ref{fig_E4D_initial} -- the lower limbs are not modeled in detail but represented by a single body that is fixed to the seat.
The joints in the upper extremities are controlled by muscles whereas the other joints are actuated using joint actuators.
In addition to this, the final hand positions on the steering wheel and the position of the upper and lower torso are constrained.
In the scenario with active backrest support, the raising of the backrest is, in addition to the movements of EMMA, part of the solution of the optimal control problem.
The weighting factors of the objective function of the optimal control problem are set to achieve a smooth, targeted and natural movement.

\begin{figure}[t]
	\centering
	\includegraphics[width=.9\linewidth]{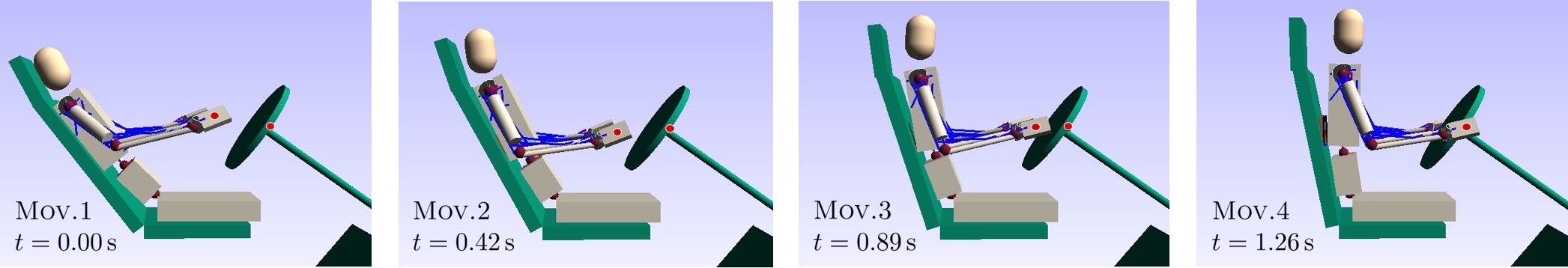}
	
	\vspace{10pt}
	
	\includegraphics[width=.9\linewidth]{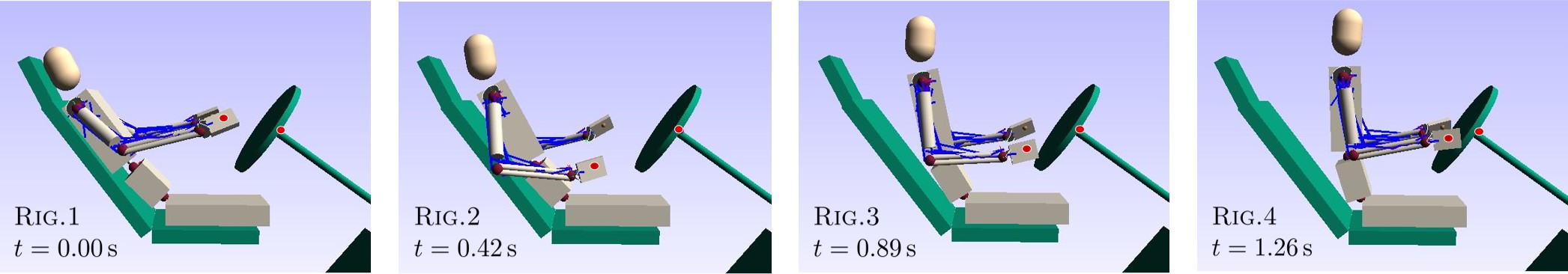}
	\caption{\enspace Reaching the steering wheel from a reclined posture. The upper frame series shows a simulation with an active backrest whereas the lower series shows a simulation with a passive backrest. In the last frame shown the passenger reaches the steering wheel only in the scenario with the support from the active backrest.\label{fig_EMMAsim}}
\end{figure}

The simulation results show that an active backrest can help to support the passenger in reaching the steering wheel in a shorter time.
The results of the simulations are obtained by solving the described optimal control problem, which takes about $45$ seconds on a standard PC and are shown in Fig.~\ref{fig_EMMAsim}.
The rows of the figure show the process of reaching the steering wheel from a reclined posture.
In the first row of the figure, the frames of the simulation with active backrest support are shown.
Hereby, the first frame shows the initial position and the last frame the final position.
It can be seen that the backrest is raised to support the passenger and the driving posture is reached within $1.3$~s.
The second row shows the frames at the same simulation times, but with passive backrest support.
Here the passenger is not able to reach the final position within the same period of time.
The steering wheel is now reached at approximately $1.6$~s.
These two simulations show that an active backrest can help reach the steering wheel faster and thus the driver can take over the driving task in a shorter time.

Although this early study is limited in its ability to realistically represent a takeover scenario in all details, it does show that active backrest support can help the driver grip the steering wheel more quickly.
While this is an important finding for the development of highly automated vehicles, it is still important to point out the special capabilities of the simulation setup.
In this study, it was possible to optimize a technical product together with the active movement of the human as a whole.
This is a unique feature of the EMMA model that cannot be achieved with motion capture methods and enables further improvement of human-machine interaction in various technical applications.

\pagebreak
\section{Discussion}\label{sec_outlook}

Vehicle safety has improved remarkably over the past century~\cite{NSC22} due, in part, to lessons learned from studying simulated crashes using mechanical and computational models of the human body.
The broad variety of design and safety considerations that are a part of modern car design have resulted in a correspondingly wide variety of digital models to represent the human form: kinematic models to assess ergonomics~\cite{BubbEtAl06,PhillipsBadler88}; multibody models to both analyze recorded motions~\cite{SethEtAl18,DamsgaardEtAl06} and simulate responses through the phases of a crash~\cite{HappeeEtAl98}; and detailed finite element models to compute the stresses and strains of ATDs~\cite{YangEtAl06} and HBMs~\cite{KatoEtAl18, CorreiaMclachlinCronin21, JohnEtAl22} throughout the crash phases.
This broad array of models and simulation methods motivated us to write this review to help people working in the field both to see what is possible and to make more informed decisions when using simulation in the future.

The three case studies we presented highlight the possibilities offered by digital models to predict in-crash phase injuries, pre-crash phase reflex responses, and voluntary movement in general. 
In our first case study, we examined the evaluation of a novel motorcycle safety system~\cite{MaierEtAl22} during the in-crash phase using a variety of ATDs and HMBs.
While the safety system dramatically reduced the impact of the head, the differences in peak head acceleration between the different types of models were larger than the differences due to anthropometry alone.
Since the risk of injury can also vary with posture in the pre-crash phase~\cite{KempterEtAl21}, our second case study focused on tuning a head-neck reflex controller specifically~\cite{KempterEtAl23} and simulating pre-crash scenarios in general.
In both cases, these tasks required multibody models, such as Madymo~\cite{HappeeEtAl03}, as the cost of performing thousands of finite element simulations would have been prohibitive.
In our final case study, we evaluated the effectiveness of an active backrest to help a driver more quickly reach the steering wheel during an emergency while autonomous driving~\cite{NHTSA13}. 
Voluntary movement is being predicted in this simulation which is strongly influenced by cognition and experience.
The final movement was arrived at by casting the task as an optimal control problem with a physiologically informed cost function~\cite{RollerEtAl17,RollerEtAl20,BjoerkenstamEtAl20} and solving the resulting optimization problem.
Since the solution process of an optimal control problem requires a large volume of simulations, the passenger was modeled with the multibody model EMMA and contact was treated using rigid constraints to make the simulations as fast as possible.
The resulting solution makes it clear that an active backrest can help the driver reach the steering wheel faster during an emergency takeover.

\begin{figure}[t]
	\centerline{\includegraphics[width=1.0\textwidth]{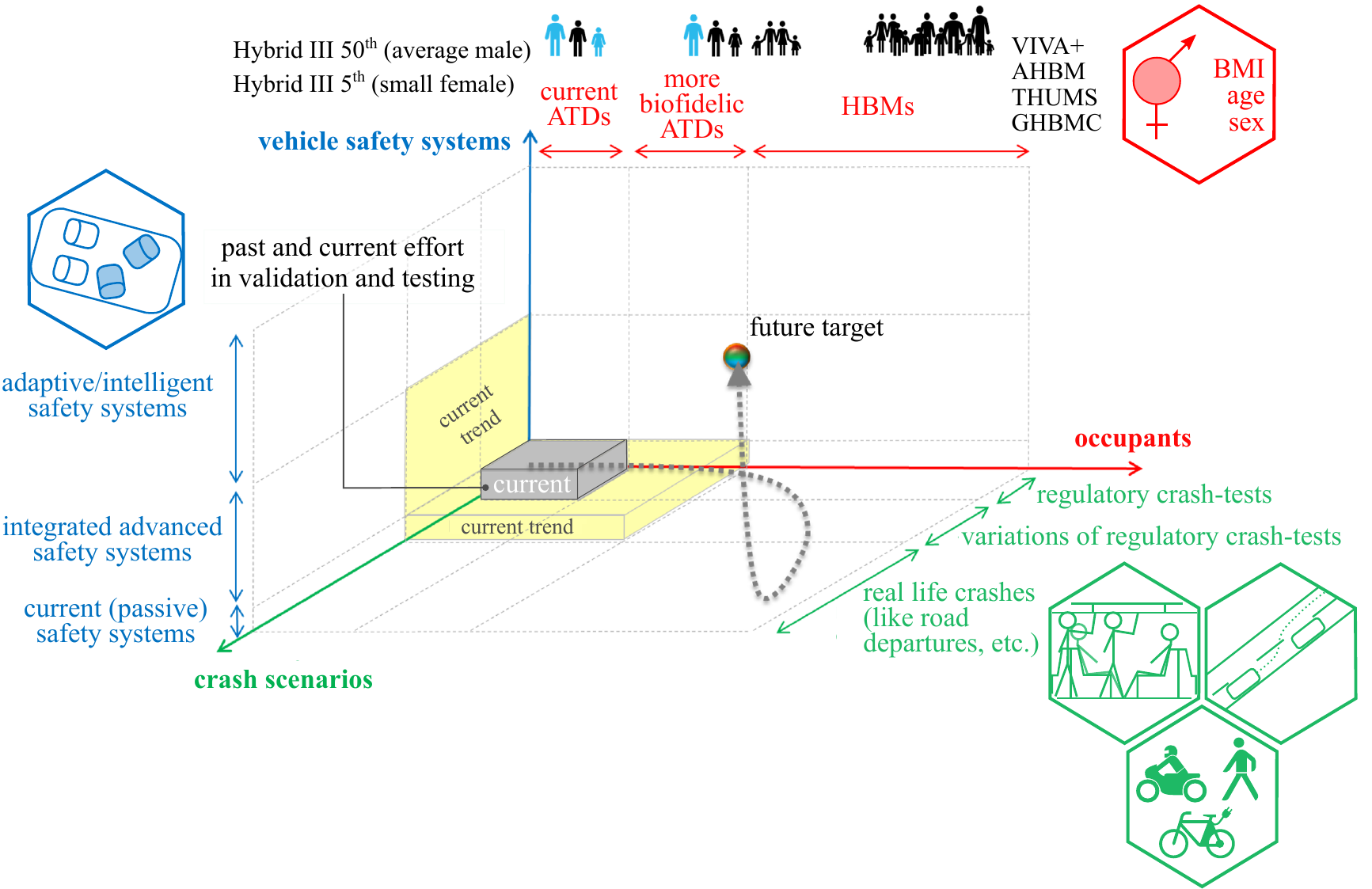}}
	\caption{\enspace Trends of the three-dimensional challenge of regulatory approval and evaluation for consumer information. The axes represent the occupant model (red), the accident scenarios (green), and the safety systems (blue) considered. The dashed line traces the path to a future target where the safety systems are evaluated considering human variability, real-life crashes, and all available safety systems.     \label{fig_ThreeDimensionalChall}}
\end{figure}

Although simulation has become an indispensable tool for designing safety systems, there are three future challenges the field of simulation will need to address (Fig.~\ref{fig_ThreeDimensionalChall}) to meet and exceed the Euro NCAP~\cite{EuroNCAP21} assessment protocol: a wider variety of models need to be developed and simulated to evaluate the risk of injury faced by a population of people (red axis in Fig.~\ref{fig_ThreeDimensionalChall}); a wider variety of crash scenarios will need to be considered as roadways become used by a greater range of electric and lightweight vehicles (green axis in Fig.~\ref{fig_ThreeDimensionalChall}); and finally, simulations will need to consider the new safety systems that autonomous vehicles introduce (blue axis in Fig.~\ref{fig_ThreeDimensionalChall}).
Most vehicle passengers are insufficiently protected~\cite{FormanEtAl19} because current safety systems have been tailored for average-sized male occupants seated in a standard position.
Industry initiatives, such as the E.V.A.~\cite{EVA} (Equal Vehicles for All), are helping to ensure that the height, weight, sex, and pregnancy status of passengers are considered during the design and evaluation of safety systems.

While HBMs are currently used to certify active bonnet safety systems for pedestrians \cite{KlugEllway21}, we expect that virtual certification will play a greater role in the future to cope with the wider variety of body types, crash scenarios, and safety systems that must be considered.
Modifications are being proposed to include virtual certification using HBMs~\cite{EuroNCAP17, EuroNCAP22} to support the increased number of tests needed to make vehicle safety equitable. 
Specialized tools are being developed~\cite{GalijatovicEtAl22} to allow a regulator to verify the integrity and authenticity of a safety simulation without compromising the confidentiality of the intellectual property contained in the models and simulations.
In addition to regulatory tools, the models used to certify a particular safety feature will likely require continued improvements to strike the right balance between detail and computational expense.

As illustrated by the three case studies, the depth of detail of a model directly influences both the calculation time and ways a model can be used: the finite element models used in Sec.~\ref{subsec:case_study_safe_motorcycle} are between 100,000 to 1,000,000 times slower than real-time, the detailed Madymo multibody model used in Sec.~\ref{subsec:case_study_safety_simulation_framework} was just over 600 times slower than real-time, while the EMMA multibody model described in Sec.~\ref{sec:case_study_EMMA4Drive} can be integrated forwards in time faster than real-time\footnote{The time to converge on a solution to the optimal control problem is only ten times slower than real-time, in part, because of EMMA's rapid forward integration time.}.
The incentives to produce accurate simulations with less computation are pushing the field to make more realistic models of human behavior, more efficient simulations of contact, and to develop surrogate models that capture the detail of a finite element model at a fraction of the simulation cost.
While it is clear that both posture and muscle activity affect the injuries sustained during a crash~\cite{KempterEtAl23}, only a few scenarios have been studied using in-vivo experiments, and it is unclear how best to incorporate these observations into models.
Leaving the complexities of the HBM aside, the run time of multibody models is often dominated by simulating contact forces between bodies.
Since passenger simulations inherently involve contact, both during typical driving and in a crash scenario, there is a strong need for contact model formulations that are both accurate and fast.
Even with a fast contact model, a finite element simulation remains the most accurate option to predict injury, yet it is computationally expensive.
However, a surrogate model can be made that can estimate the response of the finite element model at a fraction of the computational cost if a database of simulations has already been computed.

Even with more efficient models, the real-world variability attempting to be represented in simulation is daunting.
As shown in the safe motorcycle case study (Sec.~\ref{subsec:case_study_safe_motorcycle}), the results produced by a variety of models can be surprising: the peak head accelerations of the 5$^{th}$ and 50$^{th}$ percentile Hybrid~III models were quite similar, but differed dramatically from the 50$^{th}$ percentile male and female VIVA+ models.
Similarly, the differences in head acceleration were large between the Hybrid~III, VIVA+, and GHBMC 50$^{th}$ percentile male models despite having only small differences in anthropometry.
While the variations in anthropometry are often evident, it is not always clear which of the thousands of components in an HBM lead to such a variety of responses between models with similar anthropometry.
The inclusion of active muscle responses introduces an additional layer of variability that has been observed in Driver-in-the-Loop experiments, as well as in other volunteer experiments.
It remains an open question of how best to model human reflexes and validate reflex controllers to ensure the response is accurate and physiologically plausible.
Answering these questions will require the continued cooperation across the institutions and disciplines that are needed to make accurate, repeatable, and timely virtual testing using HBMs~\cite{EuroNCAP22} a reality.
Despite the challenges, virtual testing has the great potential to push the boundaries of simulation science and make vehicles safer for everyone.

\section*{Acknowledgments}
We gratefully thank Dr. Julian Hay for granting access to his results of the scenario-based \textit{pre-crash occupant simulation} as well as for numerous discussions during the preparation of this paper.

\section*{Funding Information}
The study was funded by Deutsche Forschungsgemeinschaft (DFG, German Research Foundation) under Fraunhofer and DFG-Trilateral Transfer Projects "EMMA4Drive" grant number 440904784, under Germany’s Excellence Strategy - EXC 2075 - 390740016 and EXC 310; and the state of Baden-Wuerttemberg via the Juniorprofessor program and the Innovative Mobility Solutions program. We acknowledge the support by the Stuttgart Center for Simulation Science (SimTech).

\section*{Author Contributions}

\begin{itemize}
\itemsep0em 
\item[NF]: joint main author, conceptualization, software, supervision, project administration, formal analysis, validation, visualization,  writing - original draft, writing - review, \& editing.

\item[MM]: joint main author, conceptualization, supervision, formal analysis, visualization, literature review, writing - revised version, writing - review, \& editing.

\item[SM]: writing - editing original draft, writing - review version, responsible for Case Study -- Safe Motorcycle, validation, visualization, writing - review \& editing.

\item[FK]: writing - editing original draft,  responsible for Case Study  AHBM in DIL, formal analysis, validation, visualization, conceptualization, writing - review \& editing.

\item[MR]:  writing - editing original draft, responsible for Case Study EMMA, visualization, funding acquisition, writing - review \& editing.

\item[JF]:  writing - editing original draft, formal analysis, validation, visualization, conceptualization, supervision, literature review, overview table, project administration, funding acquisition, writing - review \& editing.
\end{itemize}

\newpage
\section*{Author Biography}
\begin{biography}{\includegraphics[width=60pt]{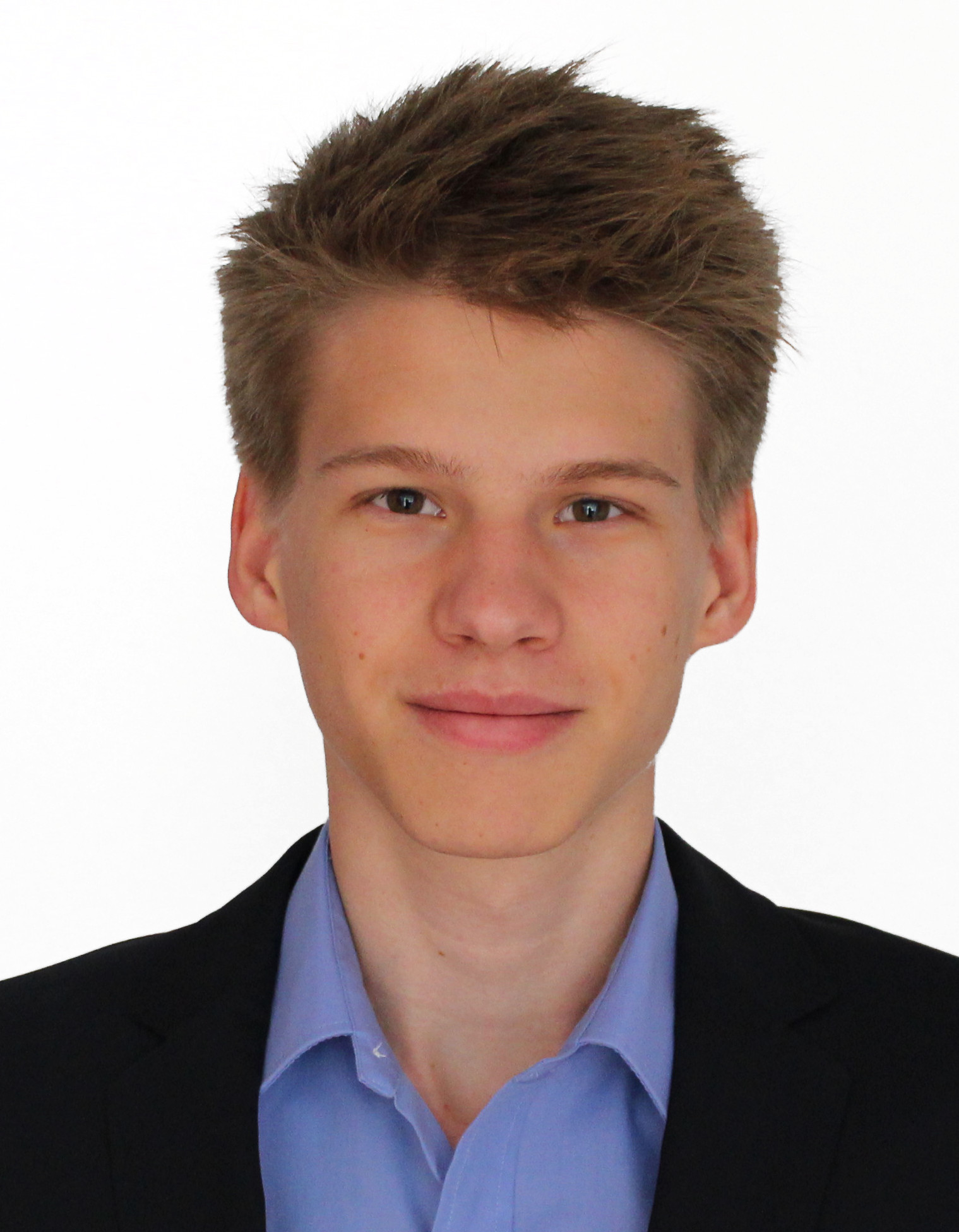}}{
	\textbf{Niklas Fahse.}
	In 2019 Niklas Fahse completed his first master's studies in Engineering Science and Mechanics at the Georgia Institute of Technology in Atlanta, Georgia, USA.
	He then proceeded with a second master's degree in Autonomous Systems at the University of Stuttgart, Germany.
	After graduating in 2021, he began doctoral studies in Mechanical Engineering at the Institute of Engineering and Computational Mechanics at the University of Stuttgart, Germany.
	His thesis will focus on surrogate models for the interaction between active human models and the vehicle interior in the context of driving simulations.
	}
\end{biography}
\begin{biography}{\includegraphics[width=60pt]{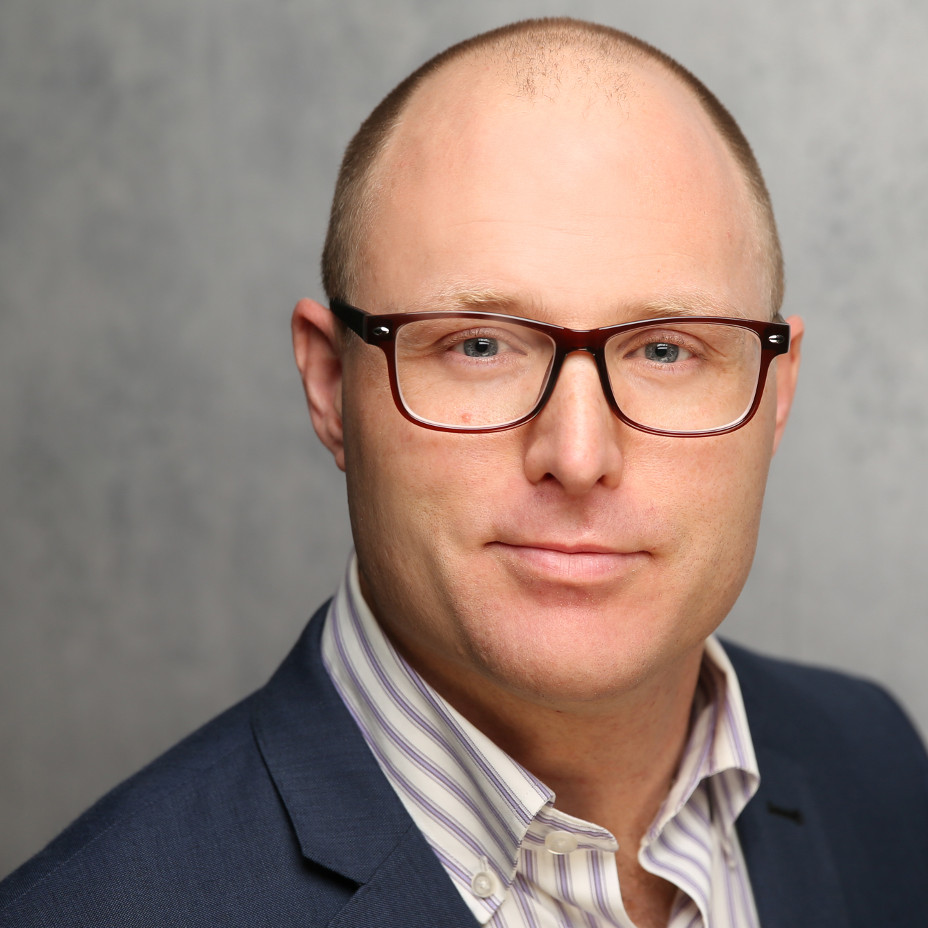}}{
		\textbf{Matthew Millard.}
		Matthew Millard completed his PhD in Systems Design Engineering in 2011 from the University of Waterloo in Canada. Following his PhD, he has taken up Postdoctoral fellowships at Stanford University, the University of Duisburg-Essen, Heidelberg University, and now works at the University of Stuttgart. During this time, he has worked on various topics, from modeling human balance to foot-ground contact modeling and predicting exoskeleton-supported lifting motions. Currently, his research focuses on the computerized prediction of neck injury during car accidents.
	}
\end{biography}
\begin{biography}{\includegraphics[width=60pt]{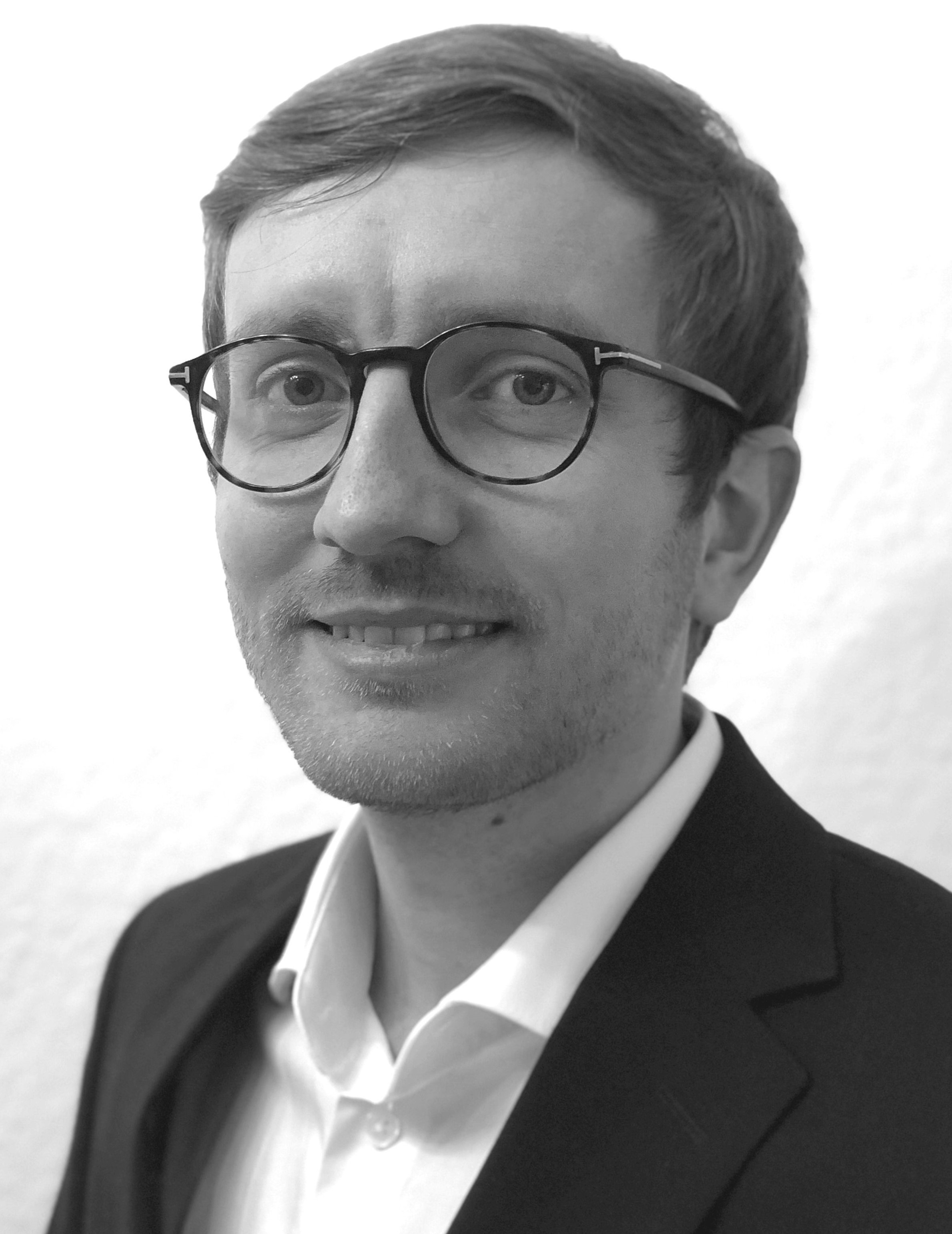}}{
	\textbf{Fabian Kempter.}
	In 2017 Fabian Kempter completed his master's degree in Mechanical Engineering at the University of Stuttgart, Germany. Subsequently, he started to pursue a PhD degree in mechanics at the Institute of Engineering and Computational Mechanics at the University of Stuttgart, Germany. The focus of his research is the study of human variability in passenger behavior and its consideration in virtual surrogate models with active musculature.
	}
\end{biography}
\vspace{0.1em}
\begin{biography}{\includegraphics[width=60pt]{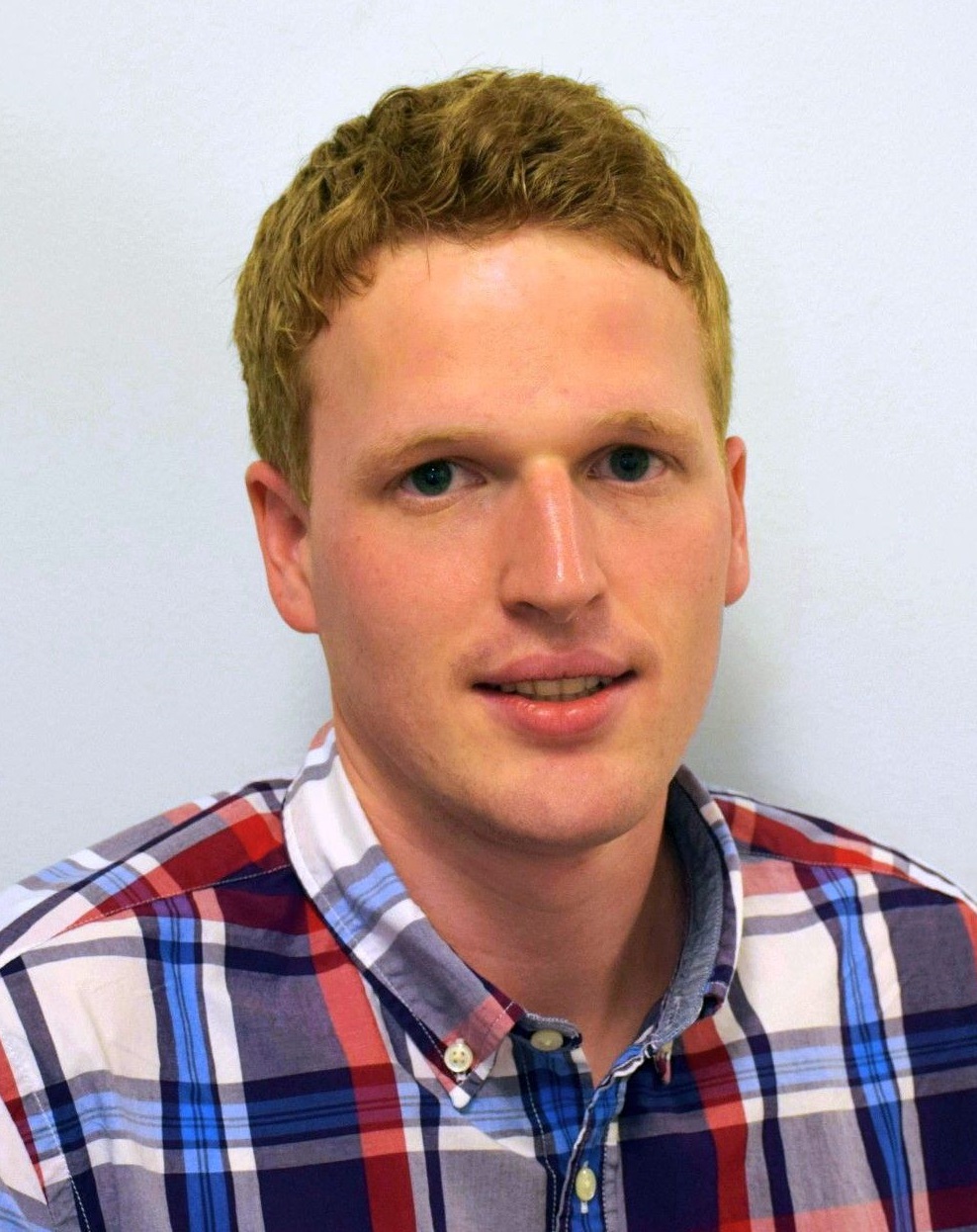}}{
	\textbf{Steffen Maier.}
	Steffen Maier is working as a research associate and doctoral student in the field of driving safety at the Institute of Engineering and Computational Mechanics (ITM) of the University of Stuttgart since 2018. He strives to improve computational modeling of powered two-wheelers and their riders as well as scientific methods to assess a rider's protection. Before, he completed his master's studies at the University of Stuttgart in Mechanical Engineering and at the Georgia Institute of Technology in Atlanta in Engineering Science and Mechanics.
	}
\end{biography}
\begin{biography}{\includegraphics[width=60pt]{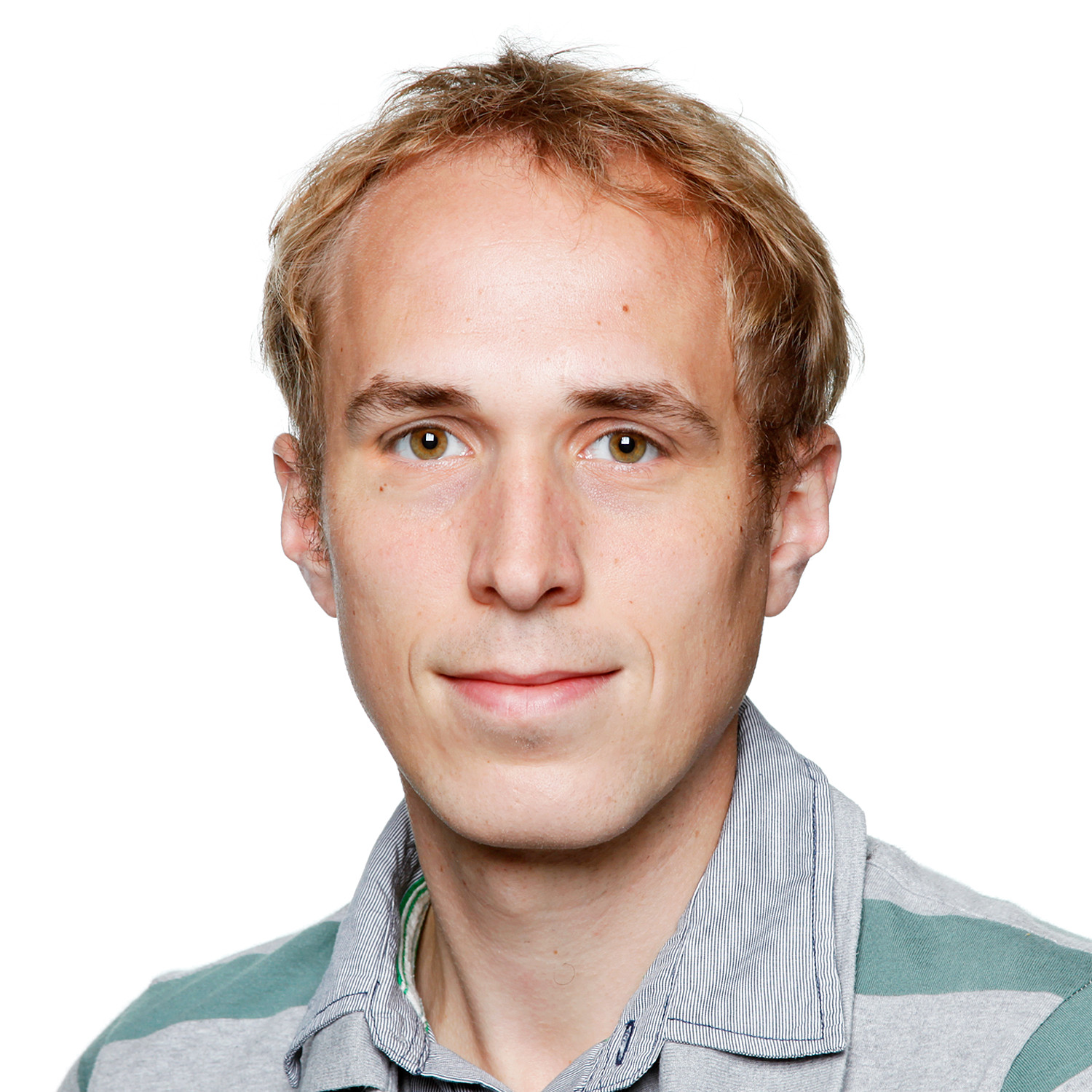}}{
	\textbf{Michael Roller.}
	Dr.-Ing Michael Roller obtained his diploma in Technomathematics at TU Kaiserslautern in 2011. He finished his doctoral degree in 2016 at Karlsruhe of Technologies in the field of structural mechanics with focus on tire simulation. From 2011 until now he is working as a researcher at the Fraunhofer Institute for Industrial Mathematics in Kaiserslautern with focus on structural mechanics and biomechanics. The focus of research in the field of digital human models is to generate physically reliable and realistic human motion without any motion capture data.
	}
\end{biography}
\begin{biography}{\includegraphics[width=60pt]{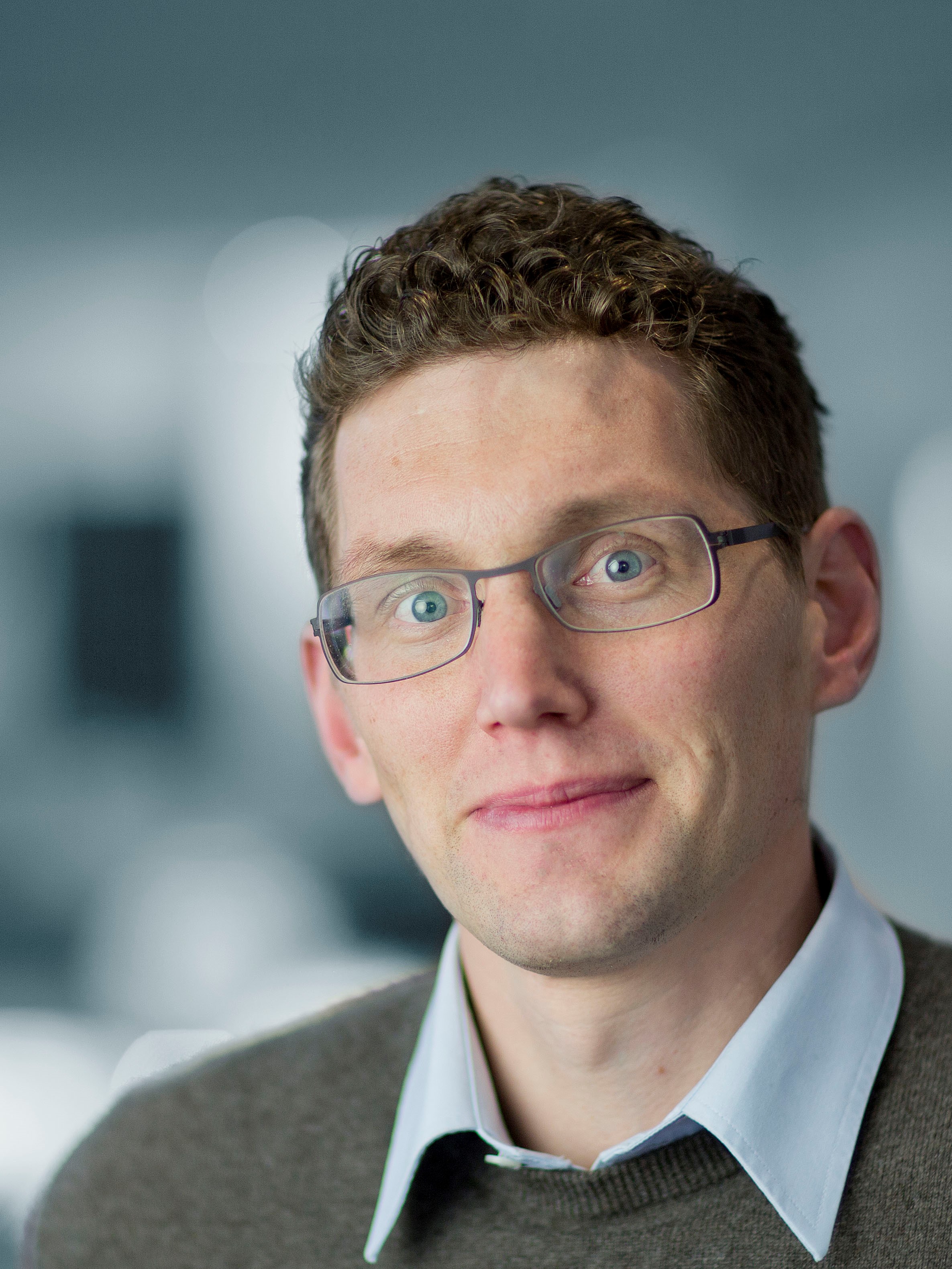}}{
	\textbf{J\"org Fehr.}
	Apl. Prof. Dr.-Ing. J\"org Fehr, studied Mechatronics in Stuttgart and Mechanical Engineering in Madison, Wisconsin, USA. He completed his PhD on Model Order Reduction in 2011 from the University of Stuttgart. After a period as a simulation engineer at TRW Automotive GmbH. In 2014 J\"org Fehr joined as Junior Professor at the Institute for Engineering and Computational Mechanics and Cluster of Excellence Simulation Technology (SRC SimTech) of the University of Stuttgart. He is on the board of Industrial Consortium SimTech. Since 2020 he is deputy director and professor after successful completion of the final evaluation. One goal of his current research is the development of optimal HBMs for the simulation of vehicle safety systems and the speedup of the simulations using linear and nonlinear model reduction methods.}
\end{biography}
\end{document}